\documentstyle[12pt,amssymb]{article}

\input epsf

\begin{document}
\baselineskip22pt

\newtheorem{prop}{Proposition}[section]
\newtheorem{theo}[prop]{Theorem}
\newtheorem{lem}[prop]{Lemma}
\newtheorem{cor}[prop]{Corollary}
\newtheorem{rk}[prop]{Remark}
\newtheorem{ex}[prop]{Example}
\newtheorem{fex}[prop]{Fundamental example}
\newtheorem{defi}[prop]{Definition}

\renewcommand {\theequation}{\arabic{section}.\arabic{equation}}

\def\pf{\noindent{\bf Proof. }}
\def\Ma#1{{\Bbb #1}}
\def\C{\Ma{C}}
\def\R{\Ma{R}}
\def\N{\Ma{N}}
\def\Z{\Ma{Z}}
\def\E{{\cal E}}
\def\F{{\cal F}}
\def\D{{\cal D}}
\def\H{{\cal H}}
\def\Gr{{\cal G}}
\def\M{{\cal M}}
\def\V{{\cal V}}
\def\I{{\cal I}}
\def\qed{\vbox{\hrule\hbox{\vrule height 1.5 ex\kern 1 ex\vrule}\hrule}}
\def\ch{\hbox{ch}\;\!}
\def\sh{\hbox{sh}\;\!}
\def\div{\hbox{div}\;\!}
\def\sgn{\hbox{sgn}}
\def\stra{{\small o}d}
\def\l{\langle }
\def\r{\rangle }
\def\p{\varphi}
\def\e{\varepsilon}
\def\o{\omega}
\def\w{\wedge}
\def\s{\sigma}
\def\a{\alpha}
\def\b{\beta}
\def\d{\delta}
\def\De{\Delta}
\def\t{\theta}
\def\g{\gamma}
\def\G{\Gamma}
\def\n{\nabla}
\def\la{\lambda}
\def\La{\Lambda}
\def\pa{\partial}
\def\til{\widetilde}
\def\iint{\int\!\!\int}
\def \non{{\nonumber}}

\title{Integration of Brownian vector fields}
\author{Yves Le Jan and Olivier Raimond} 
\date{}

\maketitle 

\begin{center}Universit\'e Paris-Sud\\ Math\'ematiques\\ B\^atiment 425\\ 91405 Orsay cedex\\ e-mail: Yves.LeJan@math.u-psud.fr; Olivier.Raimond@math.u-psud.fr\end{center}

\vskip80pt
\noindent{\bf Summary.} Using the Wiener chaos decomposition, we show that strong solutions of non Lipschitzian S.D.E.'s are given by random Markovian kernels. The example of Sobolev flows is studied in some detail, exhibiting interesting phase transitions.

\vskip60pt
\noindent{\bf Keywords~:} Stochastic differential equations, strong solution, Wiener chaos decomposition, stochastic flow, isotropic Brownian flow, coalescence, Dirichlet form.

\medskip
\noindent{\bf 1991 AMS subject classifications~:} 60H10, 31C25, 76F05.

\newpage
\setcounter{equation}{0}
\section*{Introduction.}

The purpose of this paper is to present an extended notion of strong solution to S.D.E.'s driven by Wiener processes. These solutions can be defined on rather general spaces, in the context of Dirichlet forms.

More interestingly, they are not always given by flows of maps but by flows of Markovian kernels, which means splitting can occur. Coalescent flows also appear as solutions of these S.D.E.'s. Conditions are given under which coalescence and splitting occur or not.

A variety of examples are studied. The case of isotropic Sobolev flows on the sphere or on the Euclidean space shows in particular that splitting is related to hyperinstability and coalescence to hyperstability. These notions (which will be developped in sections 9 and 10) are related to the explosion of the Lyapunov exponent toward $+\infty$ and $-\infty$. 

The typical example we have in mind is the Brownian motion on a Riemannian manifold. We consider a covariance on vector fields which induces the Riemannian metric on each tangent space. When the covariance function has enough regularity, it is known that one can solve the linear S.D.E. driven by the canonical Wiener process associated to this covariance (or in other terms to the local characteristics associated to this covariance (see section 3 below)) and get a multiplicative Brownian motion on the diffeomorphism group, which moves every point as a Brownian motion (see Le Jan and Watanabe \cite{lw} or Kunita \cite{ku}).
But models related to turbulence theory produce natural examples where the regularity condition is not satisfied. Except for the work of Darling \cite{da}, where strong solutions are not considered, these S.D.E.'s have not been really studied.
The idea is to define the solutions by their Wiener chaos expansion in terms of the heat semigroup. We call it the statistical solution. A similar expansion was given by Krylov and Veretennikov in \cite{kv}, for S.D.E.'s with strong solutions.

In this form, they appear as a semigroup of operators, and the fact that these operators are Markovian is not clearly visible in the formula. To prove this, we consider an independent realization of the Brownian motion on the manifold and couple it with the given Wiener process on vector fields using certain martingales. Then, the Markovian random operators which constitute the strong solution are obtained by filtering the Brownian motion with respect to this Wiener process. They determine the law of a canonical weak solution of the equation given the Wiener process on vector fields.
This construction has been adequately generalized to be presented in the case of symmetric diffusions on a locally compact metric space. It is a convenient and well studied framework but this assumption could clearly be relaxed (in particular to the framework of coercive forms). Relations with particle representations and filtering of S.P.D.E.'s can be observed (see Kurtz and Xiong \cite{kurtz}).

The example of Sobolev flows is studied in details on Euclidean spaces and spheres and is of major interest especially in dimension 2 and 3 where an interesting phase diagram is given in terms of the two parameters determining the Sobolev norm on vector fields : The differentiability index and the relative weight of gradients and divergence free fields (compressibility).

Some of these results have been given in the note \cite{lr} and a preliminary version of this work was released in \cite{lrb}. They are directly connected and partially motivated by a series of works of Gawedzki, Kupiainen and al on turbulent advection (\cite{bgk}, \cite{ga} and \cite{gv})

Acknowledgment~: we wish to thank an anonymous referee for his careful reading and his suggestions which helped to improve the manuscript in many ways.

\setcounter{equation}{0}
\section{Covariance function on a manifold.}

Let $X$ be a manifold.
A covariance function $C$ on $T^*X$ is a  map from $T^*X^2$ in $\R$ such that, for any $(x,y)\in X^2$, $C$ restricted to $T^*_xX \times T^*_yX$ is bilinear and such that for any $n$-uples $(\xi_1,...,\xi_n)$ of $T^*X$, 
\begin{equation}\sum_{i,j}C(\xi_i,\xi_j)\geq 0. \end{equation}

For any $\xi=(x,u)\in T^*X$, let $C_\xi$ be the vector field such that for any $\xi'=(y,v)\in T^*X$,
$$\l C_\xi(y),v\r=C(\xi,\xi').$$

Let $H_0$ be the vector space generated by the vector fields $C_\xi$. Let us define  the bilinear form on $H_0$, $\l.,.\r_H$ such that
\begin{equation} \l C_\xi,C_{\xi'}\r_H=C(\xi,\xi').\end{equation}
As equation (1.1) is satisfied, the bilinear form $\l.,.\r_H$ is a scalar product on $H_0$. We denote $\|.\|_H$ the norm associated to $\l.,.\r_H$.

Let $H$ be the separate completion of $H_0$ with respect to $\|.\|_H$. $(H,\l.,.\r_H)$ is a separable Hilbert space and we will design it as the self--reproducing space associated to the covariance function $C$. $H$ is constituted of vector fields on $X$ and for any $h\in H$ and any $\xi=(x,u)\in T^*X$,
\begin{equation} \l C_\xi,h\r_H=\l h(x),u\r.\end{equation}

Let $(e_k)_k$ be an orthonormal basis of $H$, then equation (1.3) implies that for any $\xi=(x,u)\in T^*X$,
\begin{equation} C_\xi=\sum_k\l e_k(x),u\r e_k.\end{equation}
This equation implies that for any $\xi'=(y,v)\in T^*X$,
\begin{equation} C(\xi,\xi')=\sum_k\l e_k(x),u\r\l e_k(y),v\r.\end{equation}
Therefore
\begin{equation} C=\sum_k e_k\otimes e_k.\end{equation}

\begin{rk}
On the other hand, if we start with a countable family of vector fields $(V_k)_k$, such that for any $\xi=(x,u)\in T^*X$, $\sum_k\l V_k(x),u\r^2<\infty$, it is possible to define a covariance function on $X$ by the formula
$$C=\sum_k V_k\otimes V_k.$$
\end{rk}

Examples of isotropic covariances are given in section 9 and 10. See also \cite{ba}.

Assume now a Riemannian metric $\l .,.\r_x$ is given on $X$, the linear bundles $TX$ and $T^*X$ can be identified. We will now suppose that the covariance is {\em bounded by the metric} i.e that $$C(\xi,\xi)\leq \l u,u\r_x$$ for any $\xi=(x,u)\in T^*X$. Note that this condition implies that $|h(x)|_x\leq \|h\|_H$ for any $h\in H$.

Let us denote by $m(dx)$ the volume element on $X$. Given any differentiable function $f$ such that $|\nabla f|$ is square integrable, we can map it linearly into $Df$ in the Hilbert tensor product $L^2(m)\hat{\otimes} H$ setting $\l Df,g\otimes h\rangle=\int_X g(x)\l \n f(x),h(x)\r_x~m(dx)$ for all $g\in L^2(m)$ and $h\in H$.

\smallskip
Note that 
\begin{equation}\|Df\|_H^2(x)\leq|\n f(x)|^2\end{equation}
(This notation comes from the identification $L^2(m)
\hat{\otimes} H$ with the $L^2$ space of $H$ valued functions on
$X$) and that
\begin{equation}\|Df\|^2_{L^2(m)\otimes H}\leq \int|\n f|^2dm.\end{equation}

\setcounter{equation}{0}
\section{Covariance function bounded by a Dirichlet form.}

We can extend these notions to the framework of local Dirichlet forms.
Let $X$ be a locally compact separable metric space and $m$ be a positive Radon measure on $X$ such that Supp$[m]=X$.

\bigskip
Let $(\E,\F)$ be a regular Dirichlet space, $\F\subset L^2(X,m)$. We will suppose that the Dirichlet form is local and conservative. We will denote $P_t$, the associated Markovian semigroup, $A$ its generator and $\D(A)$ the domain of $A$. We will also suppose that $m$ is an invariant measure, hence that $P_t1=1$. We will  also assume that for any $f\in\F_b=L^\infty(m)\cap \F$, there exists $\G(f,f)\in L^1(m)$ such that for any $g\in \F_b$,
\begin{equation}2\E(fg,f)-\E(f^2,g)=\int g\G(f,f)~dm.\end{equation}
$\G$ can be extended to $\F$ and we denote $\G(f,g)$ the $L^1(m)$-valued bilinear form on $\F^2$, where for any $(f,g)\in\F^2$, $\G(f,g)=\frac{1}{4}\left(\G(f+g,f+g)-\G(f-g,f-g)\right)$. A sufficient condition for the existence of $\G$ (see corollary 4.2.3 in \cite{bo}) is that $\D(A)$ contains a subspace $E$ of ${\cal D}(A)$, dense in $\F$, such that
$$\forall f\in E, \qquad f^2\in {\cal D}(A).$$
Then, for $(f,g)\in E^2$,
\begin{equation}\G(f,g)=A(fg)-f~Ag-g~Af.\end{equation}

A necessary and sufficient condition for the existence of the energy density (or carr\'e du champ operator) $\G$ is given in theorem 4.2.2 in \cite{bo}.

\begin{fex} 
$X$ is a Riemannian manifold with the metric $\l .,.\r$, $m$ is the volume measure, $\F=H^1(X)$ and for any $(f,g)\in\F^2$,
$$\E(f,g)=\frac{1}{2}\int_X \l \nabla f,\nabla g\r~dm.$$
In this case, $\G(f,g)=\l \nabla f,\nabla g\r$.
\end{fex}

Let $H$ be a separable Hilbert space and $D$ a linear map from $\F$ into the Hilbert tensor product $L^2(m)\hat{\otimes} H$ such that, for any $f\in \F$
\begin{equation}\|Df(x)\|_H^2\leq\G(f,f)(x)\end{equation}
$m(dx)$--a.e. The most interesting case is when there is equality in equation (2.3).

We define the covariance function $C$ as a bilinear map from $\F\times \F$ into $L^2(m\otimes m)$ by
\begin{equation} \l C(f,g),u\otimes v\r_{L^2(m\otimes m)}=\int_{X^2}\l Df(x),Dg(y)\r_H u(x)v(y)~m(dx)m(dy).\end{equation}
Note that 
\begin{equation} \l C(f,f),u\otimes u\r_{L^2(m\otimes m)}\leq 2\E(f,f)\|u\|^2_{L^2(m)}.\end{equation}

We will say that $C$ is a covariance function bounded by the Dirichlet form $(\E,\F)$.

\begin{rk} Alternatively, we could define the covariance as a positive bilinear map from $\F\times\F$ in $L^2(m\otimes m)$ (i.e such that for any $u_i\in L^2(m)$ and any $f_i\in\F$,
\begin{equation} \int\sum_{i,j}u_i\otimes u_j~C(f_i,f_j)~dm^{\otimes 2}\geq 0.)\end{equation}
such that
\begin{equation} \l C(f,f),u\otimes u\r_{L^2(m\otimes m)}\leq 2\E(f,f)\|u\|^2_{L^2(m)}\end{equation}
and construct as before a Hilbert space such that (2.3) holds.

Indeed, $C$ induces a linear map $\tilde{C}$ from $L^2(m)\otimes \F$ into $L^2(m)\hat{\otimes}\F$ such that
$$\l\tilde{C}(u\otimes f),v\otimes g\r_{L^2(m)\otimes \F}=\l C(f,g),u\otimes v\r_{L^2(m\otimes m)}.$$
We define $H$ as the separable closure of the space $H_0$ spanned by elements of the form $\tilde{C}(u\otimes f)$, with  $u\otimes f\in L^2(m)\otimes \F$, and equipped with the scalar product
$$\l\tilde{C}(u\otimes f),\tilde{C}(v\otimes g)\r_H=\l C(f,g),u\otimes v\r_{L^2(m\otimes m)}.$$
And for $f\in\F$, $Df$ is defined such that for any $u\otimes v\otimes g\in L^2(m)\otimes L^2(m)\otimes \F$,
$$\l Df,u\otimes \tilde{C}(v\otimes g)>_{L^2(m)\hat{\otimes}H}=\l C(f,g),u\otimes v\r_{L^2(m\otimes m)}.$$
\end{rk}

\medskip
For any $h\in H$ and $f\in \F$ define $D_hf=\l Df,h\r_H$ which belongs to $L^2(m)$. Then for any orthonormal basis $(e_k)_k$ of $H$,
\begin{equation} C=\sum_kD_{e_k}\otimes D_{e_k}. \end{equation}
Moreover, for any $f\in \F$,
\begin{equation} \|Df\|^2_H=\sum_k(D_{e_k}f)^2.\end{equation}

Remark that  condition (2.3) implies that for any finite family $(u_i,f_i)\in L^\infty(m)\times \F$,
\begin{equation}\sum_{i,j}u_iu_jD(f_i,f_j)\leq\sum_{i,j}u_iu_j\G(f_i,f_j),\end{equation}
where $D(f,g)$ denotes $\l Df(x),Dg(x)\r_H=\sum_kD_{e_k}f(x)D_{e_k}g(x)$.
When the $u_i$ are step functions with discontinuities in a set of zero measure, (2.10) is satisfied as $\sum_{i,j}u_iu_jD(f_i,f_j)=|D(\sum_iu_if_i)|^2$. Then, we can extend to any family $(u_i)$ by density in $L^2(\G(f,f)dm)$ for every $f\in\F$.

\begin{rk}It is clear that given a covariance $C$ on $T^*X$ as in Section 1,
we can build the self-reproducing space $H$ consisting of vector
fields and the mapping $D: H^1(X) \to L^2(m) \hat{\otimes} H$ so as
to construct a covariance function as in Section 2.  Now suppose
conversely that we have a separable Hilbert space $H$, a linear
map $D$ and a covariance $C$ as in Section 2, and suppose we are in
the Riemannian case.  The condition $\|Df(x)\|_H^2 \le
\Gamma(f,f)(x) = |\nabla f(x)|^2$ implies that $C(f,g)(x,y)$
depends only on $\nabla f(x)$ and $\nabla g(y)$, and so there is a
covariance $\tilde{C}$ say on $T^*X$ so that $C(f,g)(x,y) =
\tilde{C}(\nabla f(x),\nabla g(y))$. So in the Riemannian case, any
Section 2 covariance function reduces to a Section 1 covariance
function.

Further, we can now assume without any loss of generality that the
separable Hilbert space $H$ is the self-reproducing space
corresponding to $\tilde{C}$ and thus consists of vector fields.
\end{rk}

\begin{rk}
The bilinear mapping $D$ is a derivation : for any $h\in H$ and any $f\in\F$ such that $f^2\in\F$,
\begin{equation} D_hf^2=2fD_hf.\end{equation}
Note that in the Riemanian manifold case (fundamental example 2.1), $D_hf=\n_hf$ when $\G=D$.
\end{rk}

\pf We first make the remark that
$$\sum_k(D_{e_k}f^2-2fD_{e_k}f)^2=D(f^2,f^2)-4fD(f^2,f)+4f^2D(f,f).$$
Integrating this relation with respect to $m$ and using (2.10), we get that
$$\int\sum_k(D_{e_k}f^2-2fD_{e_k}f)^2~dm
\leq\int\left(\G(f^2,f^2)-4f\G(f^2,f)+4f^2\G(f,f)\right)~dm=0.$$
This implies that for every $k$, $D_{e_k}f^2-2fD_{e_k}f=0$. \qed

\setcounter{equation}{0}
\section{Construction of the statistical solutions.}

In the fundamental example 2.1, when $X$ is a Riemannian manifold, $C$ is smooth and when equality holds in (2.3), it is well known (see \cite{lw} and \cite{ku}) that a stochastic flow of diffeomorphisms on $X$ can be associated with $C$. Then, with the notations of definition 2.1 in \cite{lw}, the local characteristics of the flow are $(A,L)$ where $A=C$ and $L$ is the Laplacian on $X$.

The object of this section is to show that in the general situation considered above, it is always possible to define a flow of Markovian kernels associated with $C$ and $(\E,\F)$ (which is induced by the stochastic flow when  $C$ is smooth).

Let be given a covariance function $C$ bounded by a Dirichlet form $(\E,\F)$ on a locally compact separable metric space as in the preceding section (equation (2.3) is satisfied). Let $W_t$ be a cylindrical Brownian motion on $H$ defined on some probability space $(\Omega,{\cal A},P)$, i.e a Gaussian process indexed by $H\times \R^+$ with covariance matrix $\hbox{cov}(W_t(h),W_s(h'))=s\wedge t~\l h,h'\r_H$. Set $W_t^k=W_t(e_k)$. $(W_t^k;\;k\in\N)$ is a sequence of independent Wiener processes and  we can represent $W_t$ by $\sum_kW_t^ke_k$. Informally, the law of $W_t$ is given by
$$\frac{1}{Z}e^{-\frac{1}{2}\int_0^\infty \|\dot{W_t}\|^2_H~dt}DW.$$
Let $\F_t=\s(W_s^k;\;k\in\N;s\leq t)=\s(W_s;\;s\leq t)$.

\begin{prop}
Let $S^0_t=P_t$. We can define a sequence $S^n_t$ of random operators on $L^2(m)$ such that $E[(S^n_tf)^2]\leq P_tf^2$ in $L^1(m)$ and $S^n_t$ is ${\cal F}_t$-measurable, by the recurrence formula, in $L^2(m\otimes P)$ (i.e in the Hilbert tensor product $L^2(m)\hat{\otimes} L^2(P)$)
\begin{equation} S^{n+1}_tf=P_tf+\sum_k\int_0^tS^n_s(D_{e_k}P_{t-s}f))~dW_s^k.\end{equation}
\end{prop}

\medskip 
\noindent{\bf Remark.} The stochasic integral in equation (3.1) here makes sense as a Hilbert valued It\^o integral. Recall that given a real Wiener process $W_t$ and a Hilbert space $H$, for any $F$ progressively measurable in $L^2(P_W\otimes dt)\hat{\otimes}H$ and any $h\in H$, $\l\int_s^tF(u)dW_u,h\r_H=\int_s^t\l F(u),h\r_H dW_u$ and $E[\|\int_s^tF(u)dW_u\|_H^2]=\int_s^t\|F(u)\|_H^2du$.

\medskip

\pf Suppose we are given $S^n_t$, a $\F_t$-measurable random operator on $L^2(m)$ such that $E[(S^n_tf)^2]\leq P_tf^2$.

\smallskip
Let $f\in L^2(m)$. For any positive $t$, $P_tf\in\F$ and $D_{e_k}P_{t-s}f$ is well defined and belongs to $L^2(m)$.
\begin{eqnarray*}
E[(S^{n+1}_tf)^2] &=& (P_tf)^2+\sum_k\int_0^tE\left[\left(S^n_s(D_{e_k} P_{t-s}f)\right)^2\right]~ds\quad m-a.e.\\
&\leq& (P_tf)^2+\int_0^tP_s(|DP_{t-s}f|^2)~ds\\
&\leq& (P_tf)^2+\int_0^tP_s(\G(P_{t-s}f,P_{t-s}f))~ds.
\end{eqnarray*}
For $f\in L^{\infty}(m)\cap L^2(m)$, $\frac{\pa}{\pa s}P_s((P_{t-s}f)^2)=P_s(\G(P_{t-s}f,P_{t-s}f))$ and 
\begin{equation} P_tf^2=(P_tf)^2+\int_0^tP_s(\G(P_{t-s}f,P_{t-s}f))~ds. \end{equation}
An approximation by truncation shows that equation (3.2) remains true for $f\in L^2(m)$ and $E[(S^{n+1}_tf)^2]\leq P_tf^2$. \qed

\medskip
\noindent{\bf Remark.} {\it The definition of $S^n_t$ is independent of the choice of the basis on $H$.}

\medskip
In the following, we will use the canonical realization of the processes $W^k_t$. They will be defined as coordinate functions on $\Omega=C(\R^+,\R)^\N$, with the product Wiener measure $P$. We note $\t_t$ the natural shift on $\Omega$, such that $W^k_{t+s}-W^k_t=W^k_s\circ\t_t$.

\smallskip
Recall that an operator on $L^2(m)$ is called Markovian if and only if it preserves positivity and maps 1 into 1 (or more precisely, if $m$ is not finite, if its natural extension to positive functions maps 1 into 1). 

\begin{theo}
The family of random operators $S^n_t$ converges in $L^2(P)$ towards a one parameter family of $\F_t$--adapted Markovian operators $S_t$ such that

\smallskip
a) $S_{t+s}=S_t(S_s\circ \t_t)$, for any $s,t\geq 0$;

\smallskip
b) $\forall f\in L^2(m)$, $S_tf$ is uniformly continuous with respect to $t$ in $L^2(m\otimes P)$;

\smallskip
c) $E[(S_tf)^2]\leq P_tf^2$, for any $f\in L^2(m)$;

\smallskip
d) $S_tf=P_tf+\sum_k\int_0^tS_s(D_{e_k}P_{t-s}f)~dW_s^k$,  for any $f\in L^2(m)$;

\smallskip
e) $S_tf=f+\sum_k\int_0^tS_s(D_{e_k}f)~dW_s^k+\int_0^tS_s(Af)~ds$,  for any $f\in \D(A)$.

\smallskip
\noindent $S_t$ is uniquely characterized by c) and d) or by a), c) and e).
When $\G=D$, we call it the statistical solution of the S.D.E. (see 3.22 below)
\begin{equation}
\forall f\in\D(A):\quad df(X_t)=\sum_k D_{e_k}f(X_t)~dW_t^k+Af(X_t)~dt.
\end{equation}
Note that this S.D.E. does not always have a strong solution in the usual sense.
\end{theo}

\pf The convergence of $S^n_t$ is immediate since for any $n\geq 1$, $J^n_tf=S^n_tf-S^{n-1}_tf$ is in the Hilbert tensor product of the $n$-th Wiener chaos of $L^2(P)$ with $L^2(m)$, $S_tf=P_tf+\sum_{n=1}^\infty J^n_tf$ and $(P_tf)^2+\sum_{n\geq 1}E[(J^n_tf)^2]=\lim_{n\to\infty}E[(S^n_tf)^2]\leq P_tf^2$. It is clear that $S_t$ is $\F_t$--adapted and satisfies c). 
d) is obtained taking the limit in the recurrence formula of the proposition. 

Since
$$J^n_tf=\sum_{k_1,...,k_n}\int_{0<s_1<...<s_n<t}P_{s_1}D_{e_{k_1}}P_{s_2-s_1}\dots D_{e_{k_n}}P_{t-s_n}f~dW^{k_1}_{s_1}...dW^{k_n}_{s_n}$$
we have $J^n_{t+s}=\sum_{k\leq n}J^k_t(J^{n-k}_s\circ\t_t)$ (the $k$-th term corresponds to 
$$\sum_{k_1,...,k_n}\int_{0<s_1<...<s_k<s<s_{k+1}<...<s_n<t+s}P_{s_1}D_{e_{k_1}}P_{s_2-s_1}\dots D_{e_{k_n}}P_{t+s-s_n}f~dW^{k_1}_{s_1}...dW^{k_n}_{s_n}).$$
We deduce a) from this relation.

\smallskip
The uniqueness of a solution of d) verifying c) follows directly from the uniqueness of the Wiener chaos decomposition, obtained by iteration of d)~:
Let $T_t$ design another solution of d) and c) then for any $f\in L^2(m)$ and any integer $n$,
$$T_tf=S^{n-1}_tf+\sum_{k_1,...,k_n}\int_{0<s_1<...<s_n<t}T_{s_1}D_{e_{k_1}}P_{s_2-s_1}\dots D_{e_{k_n}}P_{t-s_n}f~dW^{k_1}_{s_1}...dW^{k_n}_{s_n}.$$
The second term of the right hand side of the preceding equation is orthogonal to the first one since its integrands are $L^2$. Indeed~:
\begin{eqnarray*}
\sum_{k_1,...,k_n}E\left[\int_{0<s_1<...<s_n<t}\left(T_{s_1}D_{e_{k_1}}P_{s_2-s_1}\dots D_{e_{k_n}}P_{t-s_n}f\right)^2~ds_1\dots ds_n\right]&\leq&\\
&&\hskip-240pt\leq\quad\sum_{k_1,\dots,k_n}\int_{0<s_1<\dots<s_n<t}P_{s_1}(|D_{e_{k_1}}P_{s_2-s_1} \dots D_{e_{k_n}}P_{t-s_n}f|^2)~ds_1\dots ds_n\\
&&\hskip-240pt\leq\quad\sum_{k_2,\dots,k_n}\int_{0<s_2<\dots<s_n<t}P_{s_2}(|D_{e_{k_2}}P_{s_3-s_2} \dots D_{e_{k_n}}P_{t-s_n}f|^2)~ds_2\dots ds_n
\end{eqnarray*}
using equation (2.3) and (3.2) and by induction is smaller than $P_tf^2$.

This proves that the Wiener chaos decomposition of $T_tf$ and $S_tf$ are the same and therefore $T_t=S_t$.

\bigskip
\noindent{\bf Proof of b).} Let us remark that for any positive $\e$, $S_{t+\e}-S_t=S_t(S_\e\circ\t_t-I)$. As $S_t$ and $S_\e\circ\t_t$ are independent and $m$ is invariant under $P_t$ , for any $f\in L^2(m)$
\begin{eqnarray}
\int E[(S_{t+\e}f-S_tf)^2]~dm&\leq& \int E[P_t(S_{\e}\circ\t_tf-f)^2]~dm\non\\
&\leq&\int E[(S_{\e}\circ\t_tf-f)^2]~dm\non\\
&\leq& \int (P_\e f^2-2fP_\e f+f^2)~dm\non\\
&\leq& 2\|f\|_{L^2(m)}\|f-P_\e f\|_{L^2(m)}.
\end{eqnarray}
Therefore, $\lim_{\e\to 0}\|S_{t+\e}f-S_tf\|_{L^2(m\otimes P)}=0$, uniformly in $t$. \qed

\begin{rk}
Note also the convergence in $L^2(m\otimes P)$ of $P_\e S_tf$ toward $S_tf$ when $\e\to 0$. Indeed $\|P_\e S_tf-S_tf\|^2_{L^2(m\otimes P)}=E[\|P_\e S_tf-S_tf\|^2_{L^2(m)}]$ and $\|P_\e S_tf-S_tf\|^2_{L^2(m)}$ converges towards 0 when $\e$ goes to 0 and is dominated by $4\|S_tf\|^2_{L^2(m)}$.
\end{rk}

\bigskip
\noindent{\bf Proof of e).} Let us remark that for any $\e$ and $t$ positive,
\begin{eqnarray}
S_{t+\e}f-S_tf &=& S_t\left(P_\e f+\sum_k\int_0^{\e}S_{u}\circ\t_t(D_{e_k}P_{\e-u}f)~dW_u^k\circ\t_t-f\right)\non\\
&=& S_t(P_\e f -f) +\sum_k\int_t^{t+\e}S_{s}(D_{e_k}P_{t+\e-s}f)~dW_s^k.
\end{eqnarray}
Hence using (3.5) for $t=\frac{i}{n}t$ and $\e=\frac{1}{n}t$, for $f\in\D(A)$, 
\begin{eqnarray*}
S_tf-f-\sum_k\int_0^tS_s(D_{e_k}f)~dW_s^k-\int_0^tS_s(Af)~ds~=\hskip-120pt&&\\
&=& \sum_{i=0}^{n-1}\left[S_{\frac{i}{n}t}(P_{\frac{t}{n}}f-f)+\sum_k\int_{\frac{i}{n}t}^{\frac{i+1}{n}t}S_{s}(D_{e_k}P_{\frac{i+1}{n}t-s}f)~dW_s^k\right.\\
&&\hskip40pt\left.-\int_{\frac{i}{n}t}^{\frac{i+1}{n}t}S_{s}(Af)~ds-\sum_k\int_{\frac{i}{n}t}^{\frac{i+1}{n}t}S_{s}(D_{e_k}f)~dW_s^k\right]\\
&=&A_1(n)+A_2(n)+A_3(n),\quad \hbox{with}
\end{eqnarray*}
\begin{eqnarray}
A_1(n)&=&\sum_{i=0}^{n-1} S_{\frac{i}{n}t}(P_{\frac{t}{n}}f-f-\frac{t}{n}Af);\\
A_2(n)&=&\sum_{i=0}^{n-1} \int_{\frac{i}{n}t}^{\frac{i+1}{n}t}(S_{\frac{i}{n}t}(Af)-S_{s}(Af))~ds;\\
A_3(n)&=& \sum_{i=0}^{n-1}\sum_k\int_{\frac{i}{n}t}^{\frac{i+1}{n}t}S_{s}(D_{e_k}(P_{\frac{i+1}{n}t-s}f-f))~dW_s^k.
\end{eqnarray}

First, using the fact that $m$ is $P_t$-invariant,
\begin{equation}\|A_1(n)\|_{L^2(m\otimes P)}\leq n\|P_{\frac{t}{n}}f-f-\frac{t}{n}Af\|_{L^2(m)}=o(1).\end{equation}

\smallskip
After, we remark that 
$$\left\|\int_{\frac{i}{n}t}^{\frac{i+1}{n}t}(S_{\frac{i}{n}t}(Af)-S_{s}(Af))~ds\right\|^2_{L^2(m\otimes P)}\leq\frac{t}{n}\int_{\frac{i}{n}t}^{\frac{i+1}{n}t}\|S_{\frac{i}{n}t}(Af)-S_{s}(Af)\|^2_{L^2(m\otimes P)}~ds.$$
As $S_t(Af)$ is uniformly continuous in $L^2(m\otimes P)$, there exists $\e(x)$ such that $\lim_{x\to 0}\e(x)=0$ and $\|S_{\frac{i}{n}t}(Af)-S_{s}(Af)\|^2_{L^2(m\otimes P)}\leq \e(\frac{t}{n})$ for any $s\in[\frac{i}{n}t,\frac{i+1}{n}t]$. Hence we get
$$\left\|\int_{\frac{i}{n}t}^{\frac{i+1}{n}t}(S_{\frac{i}{n}t}(Af)-S_{s}(Af))~ds\right\|^2_{L^2(m\otimes P)}\leq\frac{t^2}{n^2}\e(\frac{t}{n})$$
and $\|A_2(n)\|_{L^2(m\otimes P)}=o(1)$.

\smallskip
At last, as the different terms in the sum in equation (3.8) are orthogonal,
\begin{eqnarray*}
\|A_3(n)\|^2_{L^2(m\otimes P)}
&=& \sum_{i=0}^{n-1}\sum_k\int E\left[\left(\int_{\frac{i}{n}t}^{\frac{i+1}{n}t}S_{s}(D_{e_k}(P_{\frac{i+1}{n}t-s}f-f))dW_s^k\right)^2\right]dm\\
&\leq&\sum_{i=0}^{n-1}\int\int_{\frac{i}{n}t}^{\frac{i+1}{n}t}|D(P_{\frac{i+1}{n}t-s}f-f)|^2~ds~dm\\
&\leq&n\int_0^{\frac{t}{n}}\int|D(P_sf-f)|^2~dm~ds\\
&\leq&n\int_0^{\frac{t}{n}}\E(P_sf-f,P_sf-f)~ds.
\end{eqnarray*}
As $\lim_{s\to 0}\E(P_sf-f,P_sf-f)=0$, $\|A_3(n)\|_{L^2(m\otimes P)}=o(1)$.

\smallskip
Taking the limit as $n$ goes to $\infty$, this shows that $\|S_tf-f-\sum_k\int_0^tS_s(D_{e_k}f)~dW_s^k-\int_0^tS_s(Af)~ds\|_{L^2(m\otimes P)}=0$. \qed

\bigskip
\noindent{\bf Proof that a), c) and e) imply d).} Take $f\in L^2(m)$ and $\e$ positive, assuming e),
\begin{eqnarray*}
S_tP_\e f-P_tP_\e f-\sum_k\int_0^tS_s(D_{e_k}P_{t-s}P_\e f)~dW_s^k~=\hskip-210pt&&\\
&=&\sum_{i=0}^{n-1}\left[S_{\frac{i+1}{n}t}(P_{t-\frac{i+1}{n}t}P_\e f)-S_{\frac{i}{n}t}(P_{t-\frac{i}{n}t}P_\e f)-\sum_k\int_{\frac{i}{n}t}^{\frac{i+1}{n}t}S_{s}(D_{e_k}(P_{t-s}P_\e f))~dW_s^k\right]\\
&=& B_1(n)+B_2(n)+B_3(n), \quad\hbox{with}
\end{eqnarray*}
\begin{eqnarray}
B_1(n)&=&\sum_{i=0}^{n-1}\sum_k\int_{\frac{i}{n}t}^{\frac{i+1}{n}t}S_{s}(D_{e_k}(P_{t-\frac{i+1}{n}t}P_\e f-P_{t-s}P_\e f))dW_s^k;\\
B_2(n)&=&-\sum_{i=0}^{n-1}S_{\frac{i}{n}t}\left(P_{t-\frac{i}{n}t}P_\e f-P_{t-\frac{i+1}{n}t}P_\e f-\frac{t}{n}AP_{t-\frac{i+1}{n}t}P_\e f\right);\\
B_3(n)&=&\sum_{i=0}^{n-1}\int_{\frac{i}{n}t}^{\frac{i+1}{n}t}(S_{s}-S_{\frac{i}{n}t})(AP_{t-\frac{i+1}{n}t}P_\e f)~ds, \qquad\hbox{ since}\\
S_{\frac{i+1}{n}t}(P_{t-\frac{i+1}{n}t}P_\e f)&=& S_{\frac{i}{n}t}(P_{t-\frac{i+1}{n}t}P_\e f) + \int_{\frac{i}{n}t}^{\frac{i+1}{n}t}S_{s}(AP_{t-\frac{i+1}{n}t}P_\e f)~ds\non\\
&&+\sum_k\int_{\frac{i}{n}t}^{\frac{i+1}{n}t}S_sD_{e_k}P_{t-\frac{i+1}{n}t}P_\e f~dW_s^k\non.
\end{eqnarray}

\smallskip
Since the different terms in the sum in equation (3.10) are orthogonal,
\begin{eqnarray}
\|B_1(n)\|_{L^2(m\otimes P)}^2
&\leq&\sum_{i=0}^{n-1}\int\int_{\frac{i}{n}t}^{\frac{i+1}{n}t}|D(P_{t-\frac{i+1}{n}t}P_\e f-P_{t-s}P_\e f)|^2~ds~dm\non\\
&\leq&\sum_{i=0}^{n-1}\int_{\frac{i}{n}t}^{\frac{i+1}{n}t}\E(P_{t-\frac{i+1}{n}t}P_\e f-P_{t-s}P_\e f,P_{t-\frac{i+1}{n}t}P_\e f-P_{t-s}P_\e f)~ds\non\\
&\leq&n\int_0^{\frac{t}{n}}\E(P_sP_\e f-P_\e f,P_sP_\e f-P_\e f)~ds
\end{eqnarray}
as $\E(P_tf,P_tf)\leq \E(f,f)$ for any positive $t$ and any $f\in L^2(m)$. 

Equation (3.13) implies that $\|B_1(n)\|_{L^2(m\otimes P)}=o(1)$ (as $\lim_{s\to 0}\E(P_sP_\e f-P_\e f,P_sP_\e f-P_\e f)=0$).
\begin{eqnarray*}
\|B_2(n)\|_{L^2(m\otimes P)}
&\leq&\sum_{i=0}^{n-1}\|S_{\frac{i}{n}t}\left(P_{t-\frac{i}{n}t}P_\e f-P_{t-\frac{i+1}{n}t}P_\e f-\frac{t}{n}AP_{t-\frac{i+1}{n}t}P_\e f\right)\|_{L^2(m\otimes P)}\\
&\leq&\sum_{i=0}^{n-1}\left(\int\left(P_{\frac{i+1}{n}t}P_\e f-P_{\frac{i}{n}t}P_\e f-\frac{t}{n}AP_{\frac{i}{n}t}P_\e f\right)^2dm\right)^{\frac{1}{2}}\\
&\leq& n\|P_{\frac{t}{n}}P_\e f-P_\e f-\frac{t}{n}AP_\e f\|_{L^2(m)}.
\end{eqnarray*}
Hence, $\|B_2(n)\|_{L^2(m\otimes P)}=o(1)$.

\medskip
Note that if $Q_tf=E[S_tf]$, e) implies that for any $f\in\D(A)$,
$$Q_tf=f+\int_0^tQ_s(Af)ds.$$
Then $\frac{\pa}{\pa_s}Q_sP_{t-s}f=0$ for any $f\in L^2(m)$ and $0<s<t$ (then $P_{t-s}f\in\D(A))$ and we have $Q_tf=P_tf$. With this remark and the fact that a) and c) are satisfied, we see that b) and equation (3.4) are satisfied (see the proof of b). Using (3.4), we have
\begin{eqnarray*}
\|(S_s-S_{\frac{i}{n}t})(AP_{t-\frac{i+1}{n}t}P_\e f)\|^2_{L^2(m\otimes P)}~\leq\hskip-63pt&&\\
&\leq& 2\|AP_{t-\frac{i+1}{n}t}P_\e f\|_{L^2(m)}\|AP_{t-\frac{i+1}{n}t}P_\e f-P_{s-\frac{i}{n}t}AP_{t-\frac{i+1}{n}t}P_\e f\|_{L^2(m)}\\
&\leq&2\|AP_\e f\|_{L^2(m)}\|AP_\e f-P_{s-\frac{i}{n}t}AP_\e f\|_{L^2(m)}\\
&\leq&4\|AP_\e f\|^2_{L^2(m)}.
\end{eqnarray*}
Hence,
\begin{eqnarray*}
\|B_3(n)\|_{L^2(m\otimes P)}^2
&\leq&\sum_{i=0}^{n-1}\frac{t}{n}\int_{\frac{i}{n}t}^{\frac{i+1}{n}t}\|(S_s-S_{\frac{i}{n}t})(AP_{t-\frac{i+1}{n}t}P_\e f)\|^2_{L^2(m\otimes P)}~ds\\
&\leq&\frac{4t^2}{n}\|AP_\e f\|_{L^2(m)}^2.
\end{eqnarray*}
Taking the limit as $n$ goes to $\infty$, this shows that d) is satisfied for $P_\e f$, with $f\in L^2(m)$ and $\e$ positive.

\medskip
At last, since $\|S_tP_\e f-S_tf\|_{L^2(m\otimes P)}\leq \|P_\e f-f\|_{L^2(m)}$ (because c) is satisfied), $\|P_{t+\e}f-P_tf\|_{L^2(m\otimes P)}\leq \|P_\e f-f\|_{L^2(m)}$ and $\|\sum_k\int_0^tS_s(D_{e_k}P_{t-s}(P_\e f-f))dW_s^k\|_{L^2(m\otimes P)}^2\leq t\E(P_\e f-f,P_\e f-f)$. Taking the limit when $\e$ goes to $0$, we prove that d) is satisfied for any $f\in L^2(m)$. \qed

\bigskip
\noindent {\bf Proof that $S_t$ is Markovian.}

\smallskip
A more concise proof of this fact has been given in \cite{lr}, relying on Wiener exponentials and Girsanov formula. The advantage of the following proof is to be more explanatory, to give a relation with weak solutions and to yield a construction of the process law associated with the statistical solution $S_t$.

Let $(\Omega',{\cal G},{\cal G}_t,X_t,P_x)$ be a Hunt process associated to $(\E,\F)$ (see \cite{fuku}), we will take a canonical version with $\Omega'=C(\R^+,X)$. Let $W_t=\sum_kW_t^ke_k$ be a cylindrical Brownian motion on $H$, independent of the Markov process $X_t$.

Let $\M$ be the space of the martingales additive functionals, ${\cal G}_t$--adapted such that if $M\in\M$, $E_x[M_t^2]<\infty$, $E_x[M_t]=0$ q.e. and $e(M)<\infty$ where $e(M)=\sup_{t>0}\frac{1}{2t}E_m[M_t^2]$ (with $P_m=\int P_x~dm(x)$). $(\M,e)$ is a Hilbert space (see \cite{fuku}).

For $f\in\F$, $M^f\in \M$ denotes the martingale part of the semi-martingale $f(X_t)-f(X_0)$. For $g\in C_K(X)\subset L^2(\Gamma(f,f)dm)${\footnote{$C_K(X)$ design the space of functions continuous with compact support.}}, we note $g.M^f\in\M$ the martingale $\int_0^tg(X_s)~dM^f_s$, then $\M_0=\{\sum_{i=1}^ng_i.M^{f_i};\;n\in\N,\;g_i\in C_K(X),\;f_i\in\F\}$ is dense in $\M$ (see lemma 5.6.3 in \cite{fuku}), and $e(\sum_ig_i.M^{f_i})=\frac{1}{2}\sum_{i,j}\int g_ig_j\Gamma(f_i,f_j)~dm$ (see theorem 5.2.3 and 5.6.1 in \cite{fuku}).

\begin{lem}
For every $(M,N)\in\M\times\M$, there exists $\Gamma(M,N)\in L^1(m)$ such that
\begin{equation} \l M,N\r_t=\int_0^t\Gamma(M,N)(X_s)~ds, \end{equation}
where $\l.,.\r_t$ is the usual martingale bracket. And for $(f,g)\in\F$, $\Gamma(M^f,M^g)=\Gamma(f,g)$.
\end{lem}

Note that lemma 3.4 implies that $e(M,N)=\frac{1}{2}\int\G(M,N)~dm$.

\smallskip
In the fundamental example 2.1, $X_t$ is the Brownian motion on $X$, $M^f_t$ is the It\^o integral $\int_0^t\l df(X_s),dX_s\r$, $\G$ is the inverse Riemannian metric and $\M$ can be identified with the space of 1-forms equipped with the $L^2$-norm associated with the metric.

\medskip
\pf When $f\in\F$, it follows from theorem 5.2.3 in \cite{fuku} that
$$\l M^f,M^f\r_t=\int_0^t\Gamma(f,f)(X_s)~ds.$$
For $M=\sum_ih_i.M^{f_i}$, $N=\sum_jk_j.M^{g_j}$, two martingales of $\M_0$, 
\begin{equation} \l M,N\r_t=\sum_{i,j}\int_0^th_ik_j\Gamma(f_i,g_j)(X_s)~ds=\int_0^t\Gamma(M,N)(X_s)~ds,\end{equation}
with $\Gamma(M,N)=\sum_{i,j}h_ik_j\Gamma(f_i,g_j)$. $\Gamma$ is a bilinear mapping from $\M_0\times \M_0$ in $L^1(m)$. $\Gamma$ is continuous since for any $(M,N)\in\M_0\times\M_0$,
\begin{eqnarray*}
\int |\G(M,N)|~dm
&\leq& \int \G(M,M)^{\frac{1}{2}} \G(N,N)^{\frac{1}{2}}~dm\\
&\leq& 2e(M)^{\frac{1}{2}}e(N)^{\frac{1}{2}}.
\end{eqnarray*}
It follows that $\Gamma$ can be extended to $\M\times \M$.

Take $M\in\M$ and an approximating sequence $M_n\in\M_0$. Then $e(M_n-M)$ converges towards 0, $M_n$ converges towards $M$ in $L^2(P_x)$ and $\l M_n,M_n\r_t$ converges in $L^1(P_x)$ towards $\l M,M\r_t$ for almost every $x$ (see section 5-2 in \cite{fuku}). This proves that $\l M,M\r_t=\int_0^t\Gamma(M,M)(X_s)~ds$. \qed

\begin{lem}
If $m$ is bounded, for any $h\in H$, there exists a unique continuous martingale in $\M$, $N^h$ such that for any $f\in\F$, $e(N^h,M^f)=\frac{1}{2}\int D_hf~dm$ and $\frac{d}{dt}\l N^h,M^f\r_t=D_hf(X_t)$. In addition, $e(N^h)\leq \frac{1}{2}m(X)\|h\|^2$ and $\l N^h\r_t\leq \|h\|^2t$.
\end{lem}

In the Riemannian manifold case (example 2.1), $N^h_t=\int_0^t\l h(X_s),dX_s\r$ when $\G=D$.

\medskip
\pf For $h=\sum_k \la_k e_k\in H$, let us define a linear form, $\a_h$ on $\M_0$ such that for any $M=\sum_{i=1}^n g_i.M^{f_i}\in \M_0$, $\a_h(M)=\frac{1}{2}\sum_{i=1}^n\int g_iD_hf_i~dm$.
\begin{eqnarray*}
(\a_h(M))^2&=&\left(\sum_k \la_k \frac{1}{2}\sum_{i=1}^n\int g_iD_{e_k}f_i~dm\right)^2\\
&\leq& \frac{1}{4}\|h\|^2m(X)\sum_{i,j}\int g_ig_jD(f_i,f_j)~dm\leq\frac{1}{2}\|h\|^2m(X)e(M).
\end{eqnarray*}
This proves that $\a_h$ is continuous on $\M_0$ and can be extended to a continuous linear form on $\M$ such that $\a_h(M)\leq \frac{1}{\sqrt{2}}\|h\|\sqrt{m(X)e(M)}$. To this form is associated a unique $N^h\in\M$ such that $\a_h(M)=e(N^h,M)$. 

Note that for any $g\in C_K(X)$ and $f\in \F$, we have $\int gD_hf~dm=\int\Gamma(N^h,g.M^f)~dm=\int g\Gamma(N^h,M^f)~dm$. This is satisfied for every $g\in C_K(X)$, therefore for any $f\in \F$, $\Gamma(N^h,M^f)= D_hf$.

Note that we also have, for $M\in\M_0$
$$\G(N^h,M)\leq \|h\|~\G(M,M)^{\frac{1}{2}},$$
which implies that $\l N^h\r_t\leq\|h\|^2t$. \qed

\begin{rk}
When $m$ is not bounded, $N^h$ can be defined as a local martingale such that for any compact $K$ and any $f\in\F$, $1_K.N^h\in\M$, $e(1_K.N^h,M^f)=\frac{1}{2}\int_K D_hf~dm$. In addition, $e(1_K.N^h)\leq \frac{1}{2}m(K)\|h\|^2$.
\end{rk}

\medskip
Let $\gamma_{kl}$ be a function on $X$ such that $\frac{d}{dt}\l N^{e_l},N^{e_k}\r_t=\gamma_{kl}(X_t)$. Lemma 3.5 implies that the matrix $A=((\d_{kl}-\gamma_{kl}))$ is positive (as $\frac{d}{dt}\l N^h\r_t\leq\|h\|^2$). Therefore, it is possible to find a matrix $R$ such that $R^2=A$.

\begin{rk}
If for any $f\in\F$, $\Gamma(f,f)=\|Df\|^2_H$, then for any $f\in\F_b$, 
\begin{equation}M^f_t=\sum_k\int_0^tD_{e_k}f(X_s)dN_s^{e_k},\end{equation}
 $D_{e_k}f=\sum_l D_{e_l}f ~\gamma_{kl}(X_t)$ and the positive symmetric matrix $P=((\gamma_{kl}))$ is a projector. In this case, $R=I-P$.
\end{rk}

\pf Set $Q^f_t=\sum_k\int_0^tD_{e_k}f(X_s)~dN_s^{e_k}$, $Q^f\in\M$, then for any $M=\sum_{i=1}^n g_i.M^{f_i}\in \M_0$,
\begin{eqnarray*}
\l Q^f,M\r_t &=& \sum_k\int_0^tD_{e_k}f(X_s)d\l N^{e_k},M\r_s\\
&=& \sum_k\sum_{i=1}^n\int_0^tD_{e_k}f(X_s)g_i(X_s)D_{e_k}f_i(X_s)~ds\\
&=& \sum_{i=1}^n\int_0^t g_i(X_s)D(f,f_i)(X_s)~ds ~=~ \l M^f,M\r_t.
\end{eqnarray*}
This implies that for any $M\in\M$, $e(Q^f,M)=e(M^f,M)$ and $Q^f=M^f$.

\smallskip
Since by lemma 3.5, $\frac{d}{dt}\l M^f,N^{e_k}\r_t=D_{e_k}f(X_t)$, we get that
$$D_{e_k}f(X_t)=\frac{d}{dt}\l Q^f,N^{e_k}\r_t=\sum_l D_{e_l}f(X_t) \frac{d}{dt}\l N^{e_l},N^{e_k}\r_t=\left(\sum_l D_{e_l}f~\gamma_{kl}\right)(X_t).$$
This relation implies that $N^{e_k}_t=\sum_l\int_0^t\gamma_{kl}(X_s)~dN^{e_l}_s$ (this is easy to check, considering $\frac{d}{dt}\l N^{e_k},M\r_t$ with $M\in\M_0$). From this, we see that $\gamma_{kl}=\sum_i\gamma_{ki}\gamma_{il}$ (i.e $P^2=P$). \qed

\medskip
Set $\til{W}^k_t=N^{e_k}_t+\sum_l\int_0^t R_{kl}(X_s)dW^l_s$ and $\til{W}_t=\sum_k\til{W}^k_te_k$.

\medskip
In the Riemannian manifold case, when $\G(f,f)=\|Df\|_H^2$ for any $f\in\F$, denoting $C_\xi$ by $C_{(x,u)}$ when $u\in T_xX$ and $\xi=(x,u)$ we have~:
\begin{eqnarray*}
d\til{W}_t &=& dW_t+C_{(X_t,dX_t)}-C_{(X_t,dW_t(X_t))}\\
\hbox{and }\quad d\til{W}^k_t &=& dW^k_t+\l e_k(X_t),dX_t\r-\sum_l\l e_k(X_t),e_l(X_t)\r~dW^l_t.
\end{eqnarray*}
In this case, $R$ is a projector (see remark above).

\begin{lem}
$(\til{W}_t^k)_k$ is a sequence of independent Brownian motion.
\end{lem}
\pf Since $\til{W}_t^k$ is a continuous martingale, we just have to compute $\frac{d}{dt}\l \til{W}_t^k,\til{W}_t^l\r_t$ :
$$\frac{d}{dt}\l \til{W}_t^k,\til{W}_t^l\r_t=\gamma_{kl}+R^2_{kl}=\delta_{kl}.$$
This implies the lemma. \qed

\medskip
Let $\mu$ be an initial distribution of the form $hm$, with $h$ a positive function in $L^2(m)\cap L^1(m)$ and for $f\in L^2(m)$ define $\til{S}_tf$ by the conditional expectation
\begin{equation}\til{S}_tf(X_0)=E_\mu[f(X_t)|\s(X_0,\til{W}^{k}_s;\;k\in\N;\;s\leq t)].\end{equation}
(One checks easily that this definition does not depend on $h$.)
Remark that as $X_t$ is Markovian and $W_t$ has independent increments,
\begin{equation}\til{S}_tf(X_0)=E_\mu[f(X_t)|\s(X_0,\til{W}^{k}_s;\;k\in\N;\;s\geq 0)].\end{equation}

In the same way, we see that $\til{S}_t$ satisfies the multiplicative cocycle property a).

\begin{lem}
For any $f\in \D(A)$ and $\mu$ an initial distribution absolutely continuous with respect to $m$,
$$\til{S}_tf=f+\sum_k\int_0^t\til{S}_s(D_{e_k}f)~d\til{W}_s^{k}+\int_0^t\til{S}_s(Af)~ds, \quad P_\mu\hbox{ a.s}.$$
\end{lem}

\pf For any $f\in \D(A)$, we have
\begin{equation}f(X_t)=f(X_0)+M^f_t+\int_0^tAf(X_s)~ds.\end{equation}
It is clear that $E_.[\int_0^t Af(X_s)~ds|\s(\til{W}^{k}_t;\;k\in\N;\;s\leq t)]=\int_0^t\til{S}_s Af(X_s)~ds$, as (3.17) is satisfied. Let $Z_t=\sum_k\int_0^tH^k_s~d\til{W}^k_s\in L^2(\s(\til{W}^k_s;\;k\in\N;\;s\leq t))$,
\begin{eqnarray*}
E_.[Z_tM^f_t]
&=&\sum_k E_.\left[\int_0^t H^k_s ~d\l \til{W}^k,M^f\r_s\right]\\
&=&\sum_k E_.\left[\int_0^t H^k_s D_{e_k}f(X_s)~ds\right]\\
&=&\sum_k E_.\left[\int_0^t H^k_s \til{S}_s(D_{e_k}f)~ds\right]\\
&=& E_.\left[Z_t\sum_k\int_0^t \til{S}_s(D_{e_k}f)~d\til{W}_s^{k}\right].
\end{eqnarray*}
This proves that $E_.[M^f_t|\s(\til{W}^{k}_t;\;k\in\N;\;s\leq t)]=\sum_k\int_0^t \til{S}_sD_{e_k}f~d\til{W}_s^{k}$. \qed

\medskip
Now, using uniqueness in theorem 3.2 and the isomorphism $j$ between $L^2(\s(\til{W}_t^k;\;t\geq 0;\;k\in\N))$ and $L^2(\s(W_t^{k};\;t\geq 0;\;k\in\N))$, we see that $j\til{S}_t=S_t$, which implies that $S_t$ is Markovian. \qed

\begin{prop}
For any $f\in\F_b$, the martingale
\begin{equation}
P^f_t=M^f_t-\sum_k\int_0^tD_{e_k}f(X_s)d\til{W}^k_s
\end{equation}
is orthogonal to the family of martingales $\{\til{W}^k_t;\;k\in\N\}$, in the sense of the martingale bracket (i.e for any $k$, $\l P^f,\til{W}^k\r_.=0$). And for any $(f,g)\in\F_b^2$,
\begin{equation}
\l P^f,P^g\r_t=\int_0^t \left(\G(f,g)(X_s)-D(f,g)(X_s)\right)~ds.
\end{equation}
\end{prop}

\pf We just have to show that $\l P^f,\til{W}^k\r_t=0$ for every $f\in\F_b$ and every $k\in\N$ which is true as
$$\l M^f,\til{W}^k\r_t= \l M^f,N^{e_k}\r_t=\int_0^t D_{e_k}f(X_s)~ds.$$

Let $(f,g)\in\F_b^2$, then
\begin{eqnarray*}
\l P^f,P^g\r_t&=&\l P^f,M^g\r_t\\
&=&\l M^f,M^g\r_t-\sum_k\int_0^tD_{e_k}f(X_s)d\l\til{W}^k,M^g\r_s\\
&=&\int_0^t\G(f,g)(X_s)~ds-\sum_k\int_0^tD_{e_k}f(X_s)D_{e_k}g(X_s)~ds\\
&=&\int_0^t\G(f,g)(X_s)~ds-\int_0^tD(f,g)(X_s)~ds. \qquad \qed
\end{eqnarray*}

\begin{rk} In the case $\Gamma(f,f)=\|Df\|_H^2$ for any $f\in \F$, proposition 3.10 implies that $P^f_t=0$ and that
$$M^f_t=\sum_k\int_0^tD_{e_k}f(X_s)d\til{W}^k_s.$$
From this, we see that the diffusion $X_t$  satisfies the S.D.E.
\begin{equation}
f(X_t)-f(X_0) = \sum_l\int_0^t D_{e_l}f(X_s)~d\til{W}_s^k+\int_0^t Af(X_s)~ds
\end{equation}
for every $f\in\D(A)$. Therefore $(X_t,\til{W}_t)$ appears as a weak solution of this S.D.E. and $\til{S}_t$ is defined by filtering $X_t$ with respect to $\til{W}_t$.
\end{rk}

Let $P_{x,\til{\o}}(d\o')$ be the conditional law of the diffusion $X_t$, given $X_0$ and $\{\til{W}_t;\;t\in\R^+\}$ (it is independent of the choice of the initial distribution). Using the identity in law between $W$ and $\til{W}$, we get a family of conditional probabilities $P_{x,\o}(d\o')$ on $C(\R^+,X)$ defined $m\otimes P$ a.e. 

Remark that (with $X_t(\o')=\o'(t)$)
\begin{equation}
S_tf(x,\o)=\int f(X_t(\o'))~P_{x,\o}(d\o')\qquad m\otimes P \quad a.s.
\end{equation}
Under $P_{x,\o}(d\o')P(d\o)$, $X_t(\o')$ verifies the S.D.E. (3.3). It is a canonical weak solution of the S.D.E. (3.3) on a canonical extension of the probability space on which $W$ is defined. $S_t$ is obtained by filtering $X_t$ with respect to $W$.

\setcounter{equation}{0}
\section{The $n$-point motion.}

Let $P^{(n)}_t$ be the family of operators on $L^\infty(m^{\otimes n})$ such that, for any $(f_i)_{1\leq i\leq n}\in L^\infty(m)$, 
\begin{equation}P^{(n)}_t f_1\otimes ~...~ \otimes f_n=E[S_tf_1 \otimes ~...~ \otimes S_tf_n].\end{equation}
$P^{(n)}_t$ is a Markovian semigroup on $L^\infty(m^{\otimes n})$ as $S_t$ is Markovian and satisfies a) in theorem 3.2. It is easy to check that $P_t^{(2)}$ maps tensor products of $L^2(m)$ functions into $L^2(m^{\otimes 2})$.

\begin{prop}
For any family of probability laws on $X$ absolutely continuous with respect to $m$, $(\mu_i;\;1\leq i\leq n)$,
\begin{equation}
P^{(n)}_{\mu_1,...,\mu_n}(d\o'_1,...,d\o'_n)=\int_\Omega P(d\o)\otimes_{i=1}^n P_{\mu_i,\o}(d\o'_i)
\end{equation}
defines a Markov process on $X^n$ (with initial distribution $\otimes_{i=1}^m\mu_i$) associated with $P_t^{(n)}$. We shall call this Markov process on $X^n$ the $n$-point motion.
\end{prop}
\pf
For every family of functions in $L^\infty(m)$, $(f_i)_{\;1\leq i\leq n}$, $m^{\otimes n}\otimes P$ a.e (with $X^i_t(\o_i')=\o_i'(t)$)
\begin{eqnarray}
S_t^{\otimes n}f_1\otimes...\otimes f_n(x_1,...,x_n,\o)
&=& \prod_{i=1}^n S_tf_i(x_i,\o) \non \\
&=& \int\prod_{i=1}^nf_i(X^i_t(\o'_i))\otimes_{i=1}^n P_{x_i,\o}(d\o'_i).
\end{eqnarray}
We get the result by integrating both members of (4.3) with respect to $P(d\o)$. \qed

\medskip
Let $D^{(n)}$ be the linear map from $H\times \F^{\otimes n}$ in $L^2(m^{\otimes n})$ such that for any $(f_i)_{1\leq i\leq n}\in \F$ and $h\in H$,
\begin{equation}D_h^{(n)}f_1\otimes ... \otimes f_n=\sum_{i=1}^n f_1\otimes...\otimes D_hf_i\otimes...\otimes f_n.\end{equation}

\begin{prop}
For any $(f_i)_{1\leq i\leq n}\in \D(A)\cap L^\infty(m)$,
\begin{eqnarray*}
S_t^{\otimes n}f_1 \otimes...\otimes f_n&=&f_1 \otimes...\otimes f_n+\sum_k\int_0^tS_s^{\otimes n}(D_{e_k}^{(n)}f_1 \otimes...\otimes f_n)~dW_s^k\\
&&+\int_0^tS_s^{\otimes n}(A^{(n)}f_1\otimes ... \otimes f_n)~ds,
\end{eqnarray*}
where
\begin{eqnarray*}
A^{(n)}f_1\otimes ... \otimes f_n &=& \sum_{i=1}^n f_1\otimes...\otimes Af_i\otimes...\otimes f_n\\
&+&\sum_{1\leq i<j\leq n}\sum_k f_1\otimes ...\otimes D_{e_k}f_i\otimes ... \otimes D_{e_k}f_j\otimes ...\otimes f_n.
\end{eqnarray*}
\end{prop}

\noindent{\bf Remark.} 1) For $n=2$, the formula extends to functions in $\D(A)$ and $A^{(2)}f\otimes g=Af\otimes g+f\otimes Ag+C(f,g)$, where $(f,g)\in ({\cal D}(A))^2$.

2) Taking the expectation, we see that $A^{(n)}$ is the infintesimal generator of $P^{(n)}_t$ on $(\D(A)\cap L^\infty(m))^{\otimes 2}$.

3) The formula extends to $C^2_K(X^n)$ in the Riemannian manifold case (using for example the uniform density of sums of product functions and the regularizing effect of $P_\e^{\otimes n}$).

\medskip
\pf This is just a straightforward application of It\^o's formula applied to $S_tf_1 \otimes ...\otimes S_tf_n$, using the differential form of the equation satisfied by $S_t$, $e)$ in theorem 3.2. Taking the expectation and differentiating with respect to $t$, we get
\begin{eqnarray*}
\frac{d}{dt}_{\vert t=0}P^{(n)}_tf_1\otimes ... \otimes f_n
&=& \frac{d}{dt}_{\vert t=0}E[S_t^{\otimes n}f_1 \otimes...\otimes f_n]\\
&=& A^{(n)}f_1\otimes ... \otimes f_n. \qquad \qed
\end{eqnarray*}

\medskip
\begin{rk}
In general, $m^{\otimes n}$ is not invariant under $P^{(n)}_t$. 
\end{rk}

\setcounter{equation}{0}
\section{Measure preserving case.}

We say that the statistical solution $S_t$ is measure preserving if and only if $mS_t=m$ a.s for all $t$ (i.e $m$ is invariant for $S_t$).
 When $m(X)=\infty$, we use the natural extension of $S_t$ to $L^1(m)$ or to positive functions defined $m$--a.e.

Let us denote by $\F_K$ the set of functions of $\F$ which have compact support.

\begin{prop}
$S_t$ is measure preserving if and only if $\int C(f,g)~dm^{\otimes 2}$ vanishes for all $f$, $g$ in $\F_K$. Moreover, define $r_t$ on $L^2(\F_t)$ by $W_s^k\circ r_t=W^k_{t-s}-W_t^k$. Then the adjoint of $S_t$ in $L^2(m)$ is $S_t^*=S_t\circ r_t$.
\end{prop}

\begin{rk} a) When $f\in \F_K$, $C(f,f)\in L^1(m^{\otimes 2})$.

b) In the Riemannian manifold case, the condition that $\int C(f,g)~dm^{\otimes 2}$ vanishes for all $f$, $g$ in $\F_K$ is equivalent to assume that $W_t$ is divergent free in the weak sense, i.e that  for any $f\in\F_K$, $\int\l W_t,\n f\r~dm=0$. (It follows from the identity $E\left[\left(\int\l W_t,\n f\r~dm\right)^2\right]=t\int C(f,f)~dm^{\otimes 2}$.)
\end{rk}

\begin{lem}
Assume that $\int C(f,g)~dm^{\otimes 2}$ vanishes for all $f$, $g$ in $\F_K$, then for every $h\in H$, $f$, $g$ in $\F$,
\begin{equation}
\int g D_hf~dm=-\int f D_hg~dm.
\end{equation}
\end{lem}

\pf For every $h\in H$, $(g,f)\mapsto\int g D_hf~dm$ is a continuous bilinear form on $\F\times \F$ since $\|D_hf\|^2_{L^2(m)}\leq \E(f,f)\|h\|^2_{L_2(m)}$.

Take $f$, $g$ in $\F_K\cap L^\infty(m)$ then $fg\in\F_K$ (as the bounded functions of a Dirichlet space form an algebra) and, since $D_{e_k}$ is a derivation, $D_{e_k}(fg)=gD_{e_k}f+fD_{e_k}g$. Using this property, we get
\begin{eqnarray*}
\sum_k\left(\int(gD_{e_k}f+fD_{e_k}g)~dm\right)^2
&=&\sum_k\left(\int D_{e_k}(fg)~dm\right)^2\\
&=&\int C(fg,fg)~dm^{\otimes 2}=0.
\end{eqnarray*}
This implies that for every $k$, $\int g D_{e_k}f~dm=-\int f D_{e_k}g~dm$. To conclude we observe that both members of (5.1) are continuous in $f$ and $g$ and that $\F_K\cap L^\infty(m)$ is dense in $\F$ (since the Dirichlet form is regular, see section 1.1 in \cite{fuku}). \qed

\medskip
\noindent{\bf Proof of proposition 5.1.} Assume $\int C(f,g)~dm^{\otimes 2}=0$ holds for every $f$ and $g$ in $\F_K$.

Let us remark that the expression of the $n$-th chaos of $S_tf$ is given by the expression
\begin{equation}
J^n_tf=\int_{0\leq s_1\leq s_2\leq ...\leq s_n\leq t}\sum_{k_1,...,k_n} P_{s_1}D_{e_{k_1}}P_{s_2-s_1}D_{e_{k_2}}...D_{e_{k_n}}P_{t-s_n}f~dW^{k_1}_{s_1}...dW^{k_n}_{s_n}.
\end{equation}
From this expression, using lemma 5.3 and the fact that $P_t$ is self--adjoint in $L^2(m)$, we get that for $f$ and $g$ in $L^2(m)$,
\begin{eqnarray}
\int gJ^n_tf~dm &=&\\ \non
&&\hskip-55pt=~\int \int_{0\leq s_1 ...\leq s_n\leq t}f\sum_{k_1,...,k_n} (-1)^nP_{t-s_n}D_{e_{k_n}}P_{s_{n}-s_{n-1}}...P_{s_2-s_1}D_{e_{k_1}}P_{s_1}g~dW^{k_1}_{s_1}...dW^{k_n}_{s_n}~dm.
\end{eqnarray}
Making the change of variable $u_{n-i+1}=t-s_i$, we get that the adjoint of $J^n_t$ is given by
\begin{equation}
(J^n)^*_tg
=\int_{0\leq u_1\leq u_2\leq ...\leq u_n\leq t}\sum_{k_1,...,k_n} P_{u_1}D_{e_{k_1}}P_{u_2-u_1}D_{e_{k_2}}...D_{e_{k_n}}P_{t-u_n}g~dW^{k_1}_{u_1}\circ r_t...dW^{k_n}_{u_n}\circ r_t.
\end{equation}
From this it is easy to see that $S_t^*g=(S_t\circ r_t) g$ (as they have the same chaos expansion).

\smallskip
Notice that $S^*_t1=1$. A priori the constant functions are not in $L^2(m)$, but there exists an increasing sequence in $L^2(m)$, $g_n$ such that $g_n$ converges towards 1. For any nonnegative function $f\in L^2(m)$,
\begin{equation}\int S_tfg_n~dm=\int fS_t^*g_n~dm. \end{equation}
This equation implies, taking the limit as $n$ goes to $\infty$, that
\begin{equation}mS_t(f)=\int fS_t^*1~dm=m(f).\end{equation}
And we get that $mS_t=m$ a.s. Which ends the first part of the proof.

\medskip
Conversely, it follows from proposition 4.2 that for all $f$, $g$ in ${\cal D}(A)$,
$$S_t^{\otimes 2}f\otimes g-S_tf\otimes g-f\otimes S_tg+f\otimes g-\int_0^t S_s^{\otimes 2}C(f,g)~ds$$
is a square integrable martingale. This result extends to $f$, $g$ in $\F$. Taking $f$, $g$ in $\F_K$, integrating with respect to $m^{\otimes 2}$ and taking expectation, we get that $\int C(f,g)~dm^{\otimes 2}$ vanishes. \qed

\begin{rk}  
When $S_t$ is measure preserving, $P_t^{(n)}$ is self--adjoint in $L^2(m^{\otimes n})$ and in particular $m^{\otimes n}$ is invariant under $P_t^{(n)}$. The associated local Dirichlet form $\E^{(2)}$ is such that
$$\E^{(2)}(f\otimes g,f\otimes g)=\E(f,f)\|g\|^2_{L^2(m)}+\E(g,g)\|f\|_{L^2(m)}^2+2\int C(f,g)f\otimes g~dm^{\otimes 2}$$
for any $(f,g)\in\F^2$ and a similar expression can be given for $\E^{(n)}$.
\end{rk}

\setcounter{equation}{0}
\section{Existence of a flow of maps.}

Let $(S_t)_{t\geq 0}$ denote the statistical solution.

\begin{defi}
 We say that $(S_t)_{t\geq 0}$ is a flow of maps if and only if there exists a family of measurable mappings $(\p_t)_{t\geq 0}$ from $X\times \Omega$ in $X$ such that for any $f\in L^2(m)$ and any positive $t$,  $S_tf=f\circ\p_t$.
\end{defi}

Note that if $(S_t)_{t\geq 0}$ is a flow of maps, $P_{x,w}$ is the Dirac measure on the path $\{\p_t(x);\;t\geq 0\}$.

\begin{defi}
 We say that $(S_t)_{t\geq 0}$ is a coalescent flow of maps if and only if $(S_t)_{t\geq 0}$ is a flow of maps and for every $(x,y)\in X^2$, with positive probability there exists $T$ such that $\p_t(x)=\p_t(y)$ for all $t\geq T$.
\end{defi}

let $((X_t,Y_t))_{t\geq 0}$ design the two--point motion associated to the statistical solution.

\begin{defi}
We say that $(S_t)_{t\geq 0}$ is diffusive without hitting if and only if $(S_t)_{t\geq 0}$ is not a flow of maps and starting from $(x,x)$, for all positive $t$, $X_t\neq Y_t$.
\end{defi}

\begin{defi}
We say that $(S_t)_{t\geq 0}$ is diffusive with hitting if and only if $(S_t)_{t\geq 0}$ is not a flow of maps and $(X_t,Y_t)_{t\geq 0}$ hits the diagonal with positive probability.
\end{defi}

In this section, we will give conditions under which the statistical solution is a flow of maps or not.

\begin{lem} $(S_t)_{t\geq 0}$ is a flow of maps if and only if for any $f\in L^2(m)$ and any positive $t$, $E[(S_tf)^2]=P_tf^2$.
\end{lem}

\pf It is clear that there exists Markovian kernels on $X$, $s_t(x,\omega,dy)$ such that $S_tf(x)=\int f(y) s_t(x,\omega,dy)$. And $s_t(x,\o,dy)$ is the law of $X_t(\o')$ under $P_{x,\o}(d\o')$. As $m\otimes P$--a.e,
\begin{equation}(S_tf^2)(x)-(S_tf)^2(x)=\int\left(f(y)-\int f(z)s_t(x,\omega,dz)\right)^2s_t(x,\omega,dy),\end{equation}
if $E[(S_tf)^2]=P_tf^2$, $\int(f(y)-\int f(z)s_t(x,\omega,dz))^2s_t(x,\omega,dy)=0$ and $s_t(x,\omega,dz)$ is a Dirac measure $\d_{\p_t(x,\omega)}$, where $\p_t(x,\omega)$ is defined $m\otimes P$--a.e. \qed

\medskip
Let $h\in L^1(m)$ be a positive function such that $\int h~dm=1$. For any positive $t$, let $\mu_t$ be a probability on the Borel sets of $X\times X$ such that for any $(f,g)\in L^2(m)\times L^2(m)$, $\mu_t(f\otimes g)=\int E[S_tfS_tg]~h~dm$.

\begin{rk}
$(S_t)_{t\geq 0}$ is a flow of maps if and only if for all positive $t$, $\mu_t(\Delta)=1$, where $\Delta=\{(x,x);\;x\in X\}$.
\end{rk}

\pf
If $(S_t)_{t\geq 0}$ is a flow of maps, there exists $\p_t$ such that $S_tf=f\circ\p_t$. If $A$ and $B$ are disjoints Borel sets of finite measure,
$$\mu_t(A\times B)=\int E[1_A(\p_t(x))1_B(\p_t(x))]~h(x)~dm(x)=0.$$
This implies that $\mu_t(X\times X-\Delta)=0$ and as $\mu_t$ is a probability that $\mu_t(\Delta)=1$.

\smallskip
If $\mu_t(\Delta)=1$, for $f\in L^2(m)$, $\mu_t(f^2\otimes 1-2f\otimes f+1\otimes f^2)=0$. This implies that
\begin{equation} \int_X P_tf^2 ~h~dm=\int_X E[(S_tf)^2]~h~dm \end{equation}
and that $E[(S_tf)^2]= P_tf^2$. Hence $(S_t)_{t\geq 0}$ is a flow of maps. \qed
\medskip

Recall that we denoted  by $P^{(2)}_{(.,.)}$ the law of the two--point motion $((X_t,Y_t))_{t\geq 0}$.

\begin{prop}
$(S_t)_{t\geq 0}$ is a flow of maps if for any positive $r$ and any positive $t$, 
$$\lim_{y\to x}P^{(2)}_{(x,y)}[d(X_t,Y_t)\geq r]=0 \quad m(dx)-a.e.$$
\end{prop}

\pf For $\e>0$, let $\nu_{\e}$ be the measure on $X\times X$ such that for any $(f,g)\in L^2(m)\times L^2(m)$, $\nu_\e(f\otimes g)=\int f~P_\e g~h~dm$. For any $(f,g)\in L^2(m)\times L^2(m)$, 
\begin{equation}\nu_\e P^{(2)}_t(f\otimes g)=\int E[S_tfP_\e S_tg]~h~dm.\end{equation}
As $P_\e S_tg$ converges in $L^2(m\otimes P)$ towards $S_tg$ (see remark 3.3),
\begin{equation}\lim_{\e\to0}\nu_\e P^{(2)}_t(f\otimes g)=\int E[S_tfS_tg]~h~dm.\end{equation}
Therefore, the family of measure $(\nu_\e P^{(2)}_t)_{\e>0}$ converges weakly as $\e$ goes to 0 towards $\mu_t$.

\smallskip
 Assume that for any positive $r$ and any $t$, $\lim_{y\to x}P^{(2)}_{(x,y)}[d(X_t,Y_t)\geq r]=0$.  Let $A$ and $B$ be two disjoint Borel sets such that $d(A,B)\geq r$, then
$$\nu_\e P^{(2)}_t(A\times B)=\int_X f_\e(x)~h(x)~dm(x),$$
with
$$f_\e(x)=\int P^{(2)}_{(x,y)}[X_t\in A \hbox{ and } Y_t\in B]~p_\e(x,dy),$$
where $p_\e(x,dy)$ is the kernel given by $P_\e$.

\smallskip
As $d(A,B)\geq r$, 
$$f_\e(x) \leq \int P^{(2)}_{(x,y)}[d(X_t,Y_t)\geq r]~p_\e(x,dy).$$
For any positive $\b$, for $m$ almost every $x$, there exists $\a(x)$ such that $d(x,y)\leq \a(x)$ implies that $P^{(2)}_{(x,y)}[d(X_t,Y_t)\geq r]\leq\b$. Note that
$$f_\e(x) \leq \int_{\{d(x,y)>\a(x)\}}p_\e(x,dy)+\int_{\{d(x,y)\leq\a(x)\}} P^{(2)}_{(x,y)}[d(X_t,Y_t)\geq r]~p_\e(x,dy).$$

It is clear that $\lim_{\e\to 0}\int_{\{d(x,y)>\a(x)\}}p_\e(x,dy)=0$ $m(dx)$--a.e. Hence, $\limsup f_\e(x)\leq \beta$ $m(dx)$--a.e and this holds for any positive $\beta$. Therefore, $\lim_{\e\to 0}f_\e(x)=0$ $m(dx)$--a.e and by dominated convergence ($|f_\e(x)|\leq 1$) that
$$\lim_{\e\to 0}\nu_\e P^{(2)}_t(A\times B)=0.$$
This implies that $\mu_t(X\times X-\Delta)=0$ and that $(S_t)_{t\geq 0}$ is a flow of maps. \qed

\begin{prop}
If there exists a positive $t$, a positive $r$ and $p\in]0,1]$ such that for $m^{\otimes 2}$ almost every $(x,y)$, $P^{(2)}_{(x,y)}[d(X_t,Y_t)>r]\geq p$, then $(S_t)_{t\geq 0}$ is not a flow of maps.
\end{prop}

\pf Suppose there exists a positive $t$, a positive $r$ and $p\in]0,1]$ such that for $m^{\otimes 2}$ almost every $(x,y)$, $P^{(2)}_{(x,y)}[d(X_t,Y_t)>r]\geq p$.

\smallskip
Let $(B_i)_{i\in\N}$ be a partition of $X$ such that the diameter of $B_i$ is lower than $r$. 

Let us suppose that $\mu_t(\Delta)=1$ (or that $(S_t)_t$ is a flow of maps). Then we have $\sum_i\mu_t(B_i\times B_i)=1$ and for any positive $\a$, there exists $N$ such that
$$\sum_{i=1}^N\mu_t(B_i\times B_i)\geq 1-\a.$$
Since $\nu_\e P^{(2)}_t$ converges weakly towards $\mu_t$,
\begin{eqnarray*}
\sum_{i=1}^N\mu_t(B_i\times B_i)
&=& \lim_{\e\to 0}\nu_\e P^{(2)}_t(B_i\times B_i)\\
&=& \lim_{\e\to 0}\sum_{i=1}^N \int_{X\times X}P^{(2)}_{(x,y)}[(X_t,Y_t)\in B_i\times B_i]~p_\e(x,dy)~h(x)dm(x)\\
&\leq&\lim_{\e\to 0}\sum_{i=1}^N \int_{X\times X}P^{(2)}_{(x,y)}[X_t\in B_i;\;d(X_t,Y_t)\leq r]~p_\e(x,dy)~h(x)dm(x)\\
&\leq&\lim_{\e\to 0} \int_{X\times X}P^{(2)}_{(x,y)}[d(X_t,Y_t)\leq r]~p_\e(x,dy)~h(x)dm(x)\leq 1-p
\end{eqnarray*}
Choosing $\a<p$, we get a contradiction. Hence $\mu_t(\Delta)<1$ and $(S_t)_{t\geq 0}$ is not a flow of maps. \qed

\setcounter{equation}{0}
\section{A one dimensional example.}

Let $X=\R$, $P_t$ be the semigroup of the Brownian motion on $\R$ and the covariance function $C(x,y)=\sgn(x)\sgn(y)$ (where $\sgn(x)$ denotes the sign of $x$ with the convention $\sgn(0)=1$). Here, we have $W_t(x)=\sgn(x)W_t$, where $W_t$ is a Brownian motion starting from $0$. Set $L_t^{x}=\sup_{s\leq t}\{-\sgn(x)(x+W_s)\}\vee 0$ and $R_t^x=x+W_t+\sgn(x)L_t^x$ (it is a Brownian motion starting from $x$, reflected at $0$).

\begin{prop}
The statistical solution $S_t$ can be written the following way
\begin{equation}S_tf(x)=f(R_t^x)1_{L_t^x=0}+\frac{1}{2}\left[f(R_t^x)+f(-R_t^x)\right]1_{L_t^x>0}.\end{equation}
\end{prop}

\pf On an extension of the probability space, it is possible to build a Brownian motion starting from $x$, $X_t$ such that $W_t=\int_0^t\sgn(X_s)dX_s$ (then $X_t$ is a weak solution of the S.D.E. $dX_t=\sgn(X_t)dW_t$). Then $S_tf(x)=E[f(X_t)|{\cal F}^B]$, with ${\cal F}^B=\s(W_u;u\geq 0)$. Let us remark that $L_t^x$ is the local time of $X$ at $0$ and that $R_t^x=\sgn(x)|X_t|$.
Set $T=\inf\{t;~L_t^x>0\}=\inf\{t;~X_t=0\}$. The formula (7.1) follows simply from the fact that
$$E[f(X_t)1_{t\geq T}~|~|X_t|]=\frac{1}{2}\left(f(X_t)+f(-X_t)\right)1_{t\geq T}.\quad\qed$$

\setcounter{equation}{0}
\section{The Lipschitz case.}

Assume $X$ is a Riemannian manifold with injectivity radius $\rho>0$.
Let $P_t$ be the semigroup of a symmetric diffusion on $X$ with generator $A$. Let $C$ be a covariance inducing the metric (i.e with equality in (1.7)).

We will say that $C$ is Lipschitz if and only if there exist a positive constant $k$ and $0<\e<\rho$ such that~: For any $(x,y)\in X^2$, with $d(x,y)<\e$,
\begin{equation}A^{(2)}d^2(x,y)\leq k~d^2(x,y).\end{equation}

\noindent{\bf Remark.} -- $d^2(x,y)$ is smooth on $\{(x,y)\in X^2,~d(x,y)<\rho\}$ since $\rho$ is the injectivity radius.

-- On $\R^d$, the condition (8.1) will be checked as soon as $A=\frac{1}{2}\sum_{1\leq i,j\leq d}^d C^{ij}(x,x)\partial_i\partial_j+\sum_i b^i(x)\partial_i$,
\begin{equation}\sum_{i=1}^d(C^{ii}(x,x)+C^{ii}(y,y)-2C^{ii}(x,y))\leq \frac{k}{2}~d(x,y)^2\end{equation}
and $b^i$ is a Lipschitz function for all $i$. 

Equation (8.2) is  satisfied when $C$ is $C^2$ or when $C=\sum_{\a=1}^n X_\a\otimes X_\a$, where $X_\a$ are Lipschitz vector fields. In the latest case, the flow of maps can be constructed by the usual fixed point method for solutions of S.D.E.'s based on Gronwall's lemma.

\medskip
Let $(X_t,Y_t)$ be the two--point motion associated with the statistical solution. Set $\tau=\inf\{t,~d(X_t,Y_t)\geq \e\}$ and $H_t=d^2(X_{t\wedge \tau},Y_{t\wedge \tau})$.

\begin{lem}
$E^{(2)}_{(x,y)}(H_t)\leq e^{kt}~d^2(x,y)$.
\end{lem}
\pf By It\^o's formula, $$H_t-H_0=M_t+\int_0^{t\wedge\tau}A^{(2)}d^2(X_s,Y_s)~ds$$ where $M_t$ is a martingale.
Hence 
\begin{eqnarray*}
H_t-H_0&\leq& M_t+\int_0^{t\wedge\tau}k~d^2(X_s,Y_s)~ds\\
&\leq& M_t+\int_0^{t}k~H_s~ds.
\end{eqnarray*}
This implies that $E^{(2)}_{(x,y)}(H_t)-d^2(x,y)\leq k\int_0^tE^{(2)}_{(x,y)}(H_s)~ds$. Hence the lemma. \qed

\begin{theo}
Assume (8.1) is satisfied then the statistical solution associated to $P_t$ and $C$ is a flow of maps.
\end{theo}

\pf Indeed, for any $r<\e$, 
$$P^{(2)}_{(x,y)}[d(X_t,Y_t)\geq r]\leq P^{(2)}_{(x,y)}[d(X_t,Y_t)\geq r \hbox{ or }t\geq\tau]\leq \frac{1}{r^2}E^{(2)}_{(x,y)}(H_t)\leq \frac{e^{kt}}{r^2}d(x,y)^2,$$
which goes to 0 as $d(x,y)$ goes to 0. And we conclude using theorem 6.7. \qed

\setcounter{equation}{0}
\section{Isotropic statistical solution on $S^d$.}

\subsection{Isotropic covariance function on $S^d$.}

On $S^d$ with $d\geq 2$, the isotropic covariance function $C$ are given by the formula (see Raimond \cite{ra})
\begin{equation}C((x,u),(y,v))=\a(t)\l u,v\r+\b(t)\l u,y\r\l v,x\r,\end{equation}
with $(x,y)\in S^d\times S^d$, $t=\l x,y\r=\cos\p$ and $(u,v)\in T_xS^d\times T_yS^d$. $\a$ and $\b$ are given by
\begin{eqnarray}
\a(t)&=&\sum_{l=1}^\infty a_l\g_l(t)+\sum_{l=1}^\infty b_l\left(t\g_l(t)-\frac{1-t^2}{d-1}\g_l'(t)\right),\\
\b(t)&=&\sum_{l=1}^\infty a_l\g'_l(t)+\sum_{l=1}^\infty b_l\left(-\g_l(t)-\frac{t}{d-1}\g_l'(t)\right),
\end{eqnarray}
where $\g_l(t)=C_{l-1}^{\frac{d+1}{2}}(t)/C_{l-1}^{\frac{d+1}{2}}(1)$, $C_l^p$ is a Gegenbauer polynomial, $a_l$ and $b_l$ are nonnegative such that $\sum_l a_l<\infty$ and $\sum_l b_l<\infty$. Using the integral form of the Gegenbauer polynomials (see \cite{vi} p. 496):
\begin{equation}\g_l(\cos\p)=\int_0^\pi[z(\p,\t)]^{l-1}\sin^d\t~\frac{d\t}{c_d},\end{equation}
with $c_d=\int_0^\pi\sin^d\t~d\t$ and $z(\p,\t)=\cos\p-i\sin\p\cos\t$.

\medskip
In \cite{gm}, it is proved that the spectrum of the Laplacian $\De$ acting on the $L^2$--vector fields is $\{-l(l+d-1),~l\geq 1\}\cup \{-(l+1)(l+d-2),~l\geq 1\}$. Let ${\cal G}_l$ and ${\cal D}_l$ be respectively the eigenspaces corresponding to the eigenvalues $-l(l+d-1)$ and $-(l+1)(l+d-2)$. ${\cal G}_l$ is constitued of gradient vector fields and ${\cal D}_l$ of divergent free vector fields. These spaces can be isometrically identified with the spaces ${\cal H}_{d+1,l}$ and ${\cal F}_{d+1,l}$ used in \cite{ra} and can be used as carrier spaces of the irreducible representations of $SO(d+1)$, $T^l$ and $Q^l$.

Let $(\a^l_M)_M$ and $(\o^l_M)_M$ be orthonormal basis of ${\cal G}_l$ and ${\cal D}_l$. Then, if $(z^l_{M,d})_{l,M}$ and $(z^l_{M,\d})_{l,M}$ are independent families of independent normalized centered Gaussian variables,
\begin{equation}
W=\sum_{l\geq 1}\sqrt{\frac{d~a_l}{\hbox{dim}{\cal G}_l}}\sum_M z^l_{M,d}\a^l_M +\sum_{l\geq 1}\sqrt{\frac{d~b_l}{\hbox{dim}{\cal D}_l}}\sum_M z^l_{M,\d}\o^l_M
\end{equation}
is an isotropic Gaussian vector fields of covariance $C$ given by (9.1), (9.2) and (9.3).\\
{\bf Sketch of proof.} The covariance of $W$ is
$$\sum_{l\geq 1}\frac{d~a_l}{\hbox{dim}{\cal G}_l}\sum_M\a^l_M\otimes\a^l_M+\sum_{l\geq 1}\frac{d~b_l}{\hbox{dim}{\cal D}_l}\sum_M\o^l_M\otimes\o^l_M.$$
Let us choose $(\a^l_M)_M$ such that $\a^l_M=c_1(l,d)\nabla\Xi^l_M$ (where $(\Xi^l_M)_M$ is the basis of ${\cal H}_{d+1,l}$ given in \cite{ra}). Then, using the fact that $\Xi^l_M(p)=0$ if $M\neq 0$, for $x=g_1p$ and $y=g_2p$ (with $p=(0,...,0,1)$),
\begin{eqnarray*}
\sum_M\Xi^l_M(x)\Xi^l_M(y)&=&\sum_{M,N,K}T^l_{MN}(g_1)T^l_{MK}(g_2)\Xi^l_N(p)\Xi^l_K(p)\\
&=&T^l_{00}(g_2^{-1}g_1)(\Xi^l_0(p))^2.
\end{eqnarray*}
In \cite{vi} and \cite{ra}, $T^l_{00}(g)$ is computed and it is easy from this to give the covariance of the gradient part of $W$. We can calculate the covariance of the divergent free part in a similar way~: We choose the orthonormal basis $(\o^l_M)_M$ of $\D_l$ such that for $M\not\in\{1,...,d\}$, $\o^l_M(p)=0$ and such that for $1\leq i\leq d$, $\o^l_i(p)=c_2(l,d)e_i$ (this basis corresponds to the basis of $\F_{d+1,l}$ given in \cite{ra}). Then one have for $x=gp$ and $g\in SO(d)$,
\begin{equation}
\o^l_M(x) = \sum_{i=1}^d Q^l_{Mi}(g)g(\o^l_i(p))
= c_2(l,d)Q^l_{Mi}(g)g(e^i).
\end{equation}
Then for every $(x,u)$ and $(y,v)$ in $TS^d$,
\begin{eqnarray}
\sum_M \l\o^l_M(x),u\r\l\o^l_M(y),v\r
&=& (c_2(l,d))^2\sum_M Q^l_{Mi}(g_1)Q^l_{Mj}(g_2)\l g_1(e^i),u\r\l g_2(e^j),v\r\\
&=& (c_2(l,d))^2 Q^l_{ji}(g)\l g_1(e^i),u\r\l g_2(e^j),v\r,
\end{eqnarray}
with $g=g_2^{-1}g_1$. In \cite{ra}, the matrix elements $Q^l_{ji}(g)$ are calculated and it is easy from this to give the covariance of the divergence free part of $W$. \qed

\medskip
Let us now introduce Sobolev spaces and related covariances.

Let $H^{2,s}$ be the Sobolev space obtained by completion of the smooth vector fields with respect to the norm $\l (-\De +m^2)^sV,V\r_2$ (with $\l V,V\r_2=\int \|V(x)\|^2~dx$), where $m$ is positive. Note that the definition of $H^{2,s}$ does not depend on $m$. 

Let $a$ and $b$ be nonnegative reals. Take $a_l=\frac{a}{(l-1)^{\a+1}}$ and $b_l=\frac{b}{(l-1)^{\a+1}}$ for $l\geq 1$ and $a_1=b_1=0$. 
For $\a>0$, set $G(\p)=\sum_{l\geq 2}\frac{1}{(l-1)^{\a+1}}\g_l(\cos\p)$. The function $G$ is well defined on $[0,\pi]$ as $|\g_l|\leq 1$.

Let $F_d$ and $F_\d$ be real functions such that for all $l\geq 2$
\begin{eqnarray}
(l-1)^{\a+1}\dim\Gr_l\times F_d(-l(l+d-1))&=&d\\
(l-1)^{\a+1}\dim\D_l\times F_\d(-(l+1)(l+d-2))&=&d
\end{eqnarray}
and $F_d(-d)=F_\d(-2(d-1))=0$. Note that when $d=2$, $F_d=F_\d$.

Let $\Pi$ be the orthonormal projection on the space of the $L^2$--gradient vector fields.

\begin{prop}The covariance function defined by the sequences $(a_l)$ and $(b_l)$ is given by (9.1) with the functions
\begin{eqnarray}
\a(\cos\p)&=&aG(\p)+b\left(\cos\p~G(\p)+\frac{\sin\p}{d-1}\times G'(\p)\right),\\
\b(\cos\p)&=&-\frac{a}{\sin\p}G'(\p)+b\left(-G(\p)+\frac{\cos\p}{(d-1)\sin\p}\times G'(\p)\right).
\end{eqnarray}
When $a$ and $b$ are positive, the associated self--reproducing space is $H^{2,\frac{\a+d}{2}}$ equipped with a different (but equivalent) norm, namely
$$\|V\|_H^2=\frac{1}{a}\|\Pi V\|^2_d+\frac{1}{b}\|(I-\Pi)V\|_\d^2,$$where $\|V\|^2_d=\l F_d(\De)^{-1}V,V\r_2$ and $\|V\|^2_\d=\l F_\d(\De)^{-1}V,V\r_2$.
\end{prop}

\medskip
\pf It is not difficult to see that the norm $\|.\|_H$ given in the proposition is the norm on the self-reproducing space associated to $C$. 

Now since (see \cite{gm})
\begin{eqnarray*}
\dim\Gr_l&=&\frac{(d+l-3)!}{(d-1)!(l-1)!}(d+2l-3)(d+1),\\
\dim\D_l&=&\frac{(d+l-3)!}{(d-1)!(l-1)!}(d+2l-3)\frac{d(d+1)}{2},
\end{eqnarray*}
for $\la\to\infty$, $\la^{\frac{\a+d}{2}}F_d(\la)=O(1)$ and $\la^{\frac{\a+d}{2}}F_\d(\la)=O(1)$. This implies that $\|.\|_H$ and the norm used to define $H^{2,\frac{\a+d}{2}}$ are equivalent (when $a$ and $b$ are positive). And we get that the self--reproducing space associated to $C$ is $H^{2,\frac{\a+d}{2}}$. \qed 

\begin{rk} If $a$ or $b$ vanishes, the self--reproducing space is $H^{2,\frac{\a+d}{2}}$ restricted to divergent free vector fields or gradient vector fields.
\end{rk}

\subsection{Phase transitions for the Sobolev statistical solution.}
Let $P_t$ be the semigroup of the Brownian motion of variance $(a+b)G(0)$ and $S_t$ be the statistical solution associated to $P_t$ and $C$.

Let $(X_t,Y_t)$ be the two--point motion. Let $\psi_t=d(X_t,Y_t)$.
Since $h(x,y)=d^2(x,y)$ is a $C^2$-function, $h$ belongs to $\D(A^{(2)})$ and since $X_t$ and $Y_t$ are solutions of an S.D.E. like (3.3), $\psi_t^2$ is a diffusion on $[0,\pi^2]$ and is solution of an S.D.E. $\psi_t$ is also a diffusion on $[0,\pi]$ (note that $d(x,y)$ a priori does not belong to $\D(A^{(2)})$). This diffusion is eventually reflected (or absorbed) in 0 and $\pi$. Its generator is $L=\s^2(\p)\frac{d^2}{d\p^2}+b(\p)\frac{d}{d\p}$ (see Raimond \cite{ra}), with
\begin{eqnarray}
\s^2(\p)&=&\a(1)-\a(\cos\p)\cos\p+\b(\cos\p)\sin^2\p,\\
b(\p)&=&\frac{(d-1)}{\sin\p}(\a(1)\cos\p-\a(\cos\p)).
\end{eqnarray}
The generator of $\psi_t^2$ is $L'=\til{\s}^2(x)\frac{d^2}{dx^2}+\til{b}(x)\frac{d}{dx}$, with
\begin{eqnarray}
\til{\s}^2(x)&=&4x\s^2(\sqrt{x}),\\
\til{b}(x)&=&2\s^2(\sqrt{x})+2\sqrt{x}b(\sqrt{x}).
\end{eqnarray}

\begin{lem}
If $\a> 2$,  the statistical solution is a flow of maps. 
\end{lem}

\pf We have $A^{(2)}d^2(x,y)=2\s^2(d(x,y))+2b(d(x,y))d(x,y)$. When $\a>2$, then $G$ is $C^2$, this implies that $\a$ is $C^2$ and $\b$ is continuous. Hence equation (8.1) can be checked. \qed

\medskip
Suppose $a+b>0$ and let $\eta=\frac{b}{a+b}$.

\begin{theo} For any $\a\in]0,2[$,
\begin{itemize}
\item For $d=2$ or $3$ and $\eta<1-\frac{d}{\a^2}$, the statistical solution is a coalescent flow of maps.
\item For $d=2$ or $3$ and $1-\frac{d}{\a^2}<\eta<\frac{1}{2}-\frac{(d-2)}{2\a}$, the statistical solution is diffusive with hitting.
\item For $d=2$ or $3$ and $\eta>\frac{1}{2}-\frac{(d-2)}{2\a}$ or for $d\geq 4$, the statistical solution is diffusive without hitting.
\end{itemize}
\end{theo}

\noindent{\bf Remark.} The same phase transition appears in the $\R^d$ case (see theorem 10.1 below). It has been independently observed, in the context of the advection of a passive scalar, by Gawedzky and Vergassola \cite{gv}.

\begin{lem} For $\a\in ]0,2[$, we have
\begin{itemize} 
\item $G$ is differentiable on $]0,\pi[$.
\item $\lim_{\p\to 0+}\frac{G(0)-G(\p)}{\p^\a}=\int_0^\pi\int_0^\infty\frac{\cos^2\t}{t^2+\cos^2\t}t^{\a-1}\sin^d\t~\frac{dt~d\t}{\G(\a+1)c_d}=KG(0).$
\item $\lim_{\p\to 0+}\frac{G'(\p)}{\p^{\a-1}}=-\a\int_0^\pi\int_0^\infty\frac{\cos^2\t}{t^2+\cos^2\t}t^{\a-1}\sin^d\t~\frac{dt~d\t}{\G(\a+1)c_d}=-\a KG(0).$
\end{itemize}
\end{lem}

The proof of lemma 9.5 is in appendix A. From this lemma, we get as $\p$ goes to $0$
\begin{eqnarray}
\a(\cos\p)&=&(a+b)G(0)-\left(a+\left(1+\frac{\a}{d-1}\right)b\right)KG(0)\p^\a+o(\p^\a)\\
\b(\cos\p)&=&\a\left(a-\frac{b}{d-1}\right)KG(0)\p^{\a-2}+o(\p^{\a-2}).
\end{eqnarray}
Hence, 
\begin{eqnarray}
\s^2(\p)&=&(a+b)KG(0)(\a+1-\a\eta)\p^\a(1+o(1))\\
b(\p)&=&(a+b)KG(0)(d-1+\a\eta)\p^{\a-1}(1+o(1)).
\end{eqnarray}

\medskip
In order to prove theorem 9.4, we need to study the two--point motion. Because of isotropy, it is enough to study the diffusion $\psi_t$. This diffusion satisfies an S.D.E. until it exits $]0,\pi[$.

Let $s$ be the scale function of the diffusion $\psi_t$,
$$s(x)=\int_{x_0}^x\exp\left[-\int_{x_0}^y\frac{b(\p)}{\s^2(\p)}~d\p\right]~dy, \hbox{ with } (x_0,x)\in ]0,\pi[^2.$$

Let $x\in\{0,\pi\}$ and $T_x=\inf\{t>0;\;\psi_t=x\}$. Using Breiman's terminology (see \cite{br} p.368-369), $x$ is an open boundary point if $T_x=\infty$ and is a closed boundary point if $T_x<\infty$. Note that $x$ is an open boundary point if $|s(x)|=\infty$ .

\medskip
Firstly we are going to show that $\pi$ is an open boundary point. Then~:
\begin{itemize}
\item When $d=2$ or 3 and $\eta<1-\frac{d}{\a^2}$, we prove that 0 is an exit boundary point (this implies that the statistical solution is a coalescent flow of maps).
\item When $d=2$ or 3 and $1-\frac{d}{\a^2}<\eta<\frac{1}{2}-\frac{(d-2)}{2\alpha}$, we prove that 0 is an instantaneously reflecting regular boundary point (this implies that the statistical solution is diffusive with hitting).
\item When $\eta>\frac{1}{2}-\frac{(d-2)}{2\alpha}$, we prove that 0 is an open entrance boundary point (this implies that the statistical solution is diffusive without hitting).
\end{itemize}

\begin{lem} $\pi$ is an open boundary point.
\end{lem}
\pf It is easy to check that $s(\pi-)=\infty$ using the fact that $\a(1)+\a(-1)> 0$~:
$$\a(1)+\a(-1)=(a+b)G(0)+(a-b)G(\pi)>(a+b)G(\pi)+(a-b)G(\pi)\geq 0. \quad\qed$$

Since $\pi$ is an open boundary point, we now study the behavior of $\psi_t$ at and near 0.

\begin{lem} If $\eta>\frac{1}{2}-\frac{(d-2)}{2\a}$, $s(0+)=-\infty$ and if $\eta<\frac{1}{2}-\frac{(d-2)}{2\a}$, $s(0+)>-\infty$.\end{lem}
\pf Let us note $\mu=\frac{d-1+\a\eta}{\a+1-\a\eta}$. Then, we have that $\frac{b(\p)}{\s^2(\p)}=\frac{\mu}{\p}(1+o(1))$ and for any positive $\e$ there exist positive constants $C_1$ and $C_2$ such that for $y\leq x_0$,
\begin{equation}
C_1y^{-\mu+\e}\leq\exp\left[-\int_{x_0}^y\frac{b(\p)}{\s^2(\p)}~d\p\right]\leq C_2 y^{-\mu-\e}.
\end{equation}
From this, we see that $s(0+)=-\infty$ if $\mu>1$ (or if $\eta>\frac{1}{2}-\frac{(d-2)}{2\a}$) and $s(0+)$ is finite if $\mu<1$ (or if $\eta<\frac{1}{2}-\frac{(d-2)}{2\a}$). \qed

\medskip
Lemma 9.6 and 9.7 implies that (see theorem VI-3.1 in \cite{iw}) if $\eta<\frac{1}{2}-\frac{(d-2)}{2\a}$ we have $T_0<\infty$, $T_\pi=\infty$ a.s. and if $\eta>\frac{1}{2}-\frac{(d-2)}{2\a}$, 0 is an open boundary point and we have $\liminf\psi_t=0$ and $\limsup\psi_t=\pi$ a.s ($\psi_t$ is recurrent).
\begin{rk} When $d\geq 4$ and $\a\in]0,2[$, $\frac{1}{2}-\frac{(d-2)}{2\a}<0$. This implies that $\liminf\psi_t=0$ and $\limsup\psi_t=\pi$ a.s.\end{rk}

Since $\pi$ is an open boundary point, $\psi_t\in [0,\pi[$ for every positive $t$ and $\psi^2_t$ is a solution of the S.D.E.
\begin{equation}
d\psi_t^2=\sqrt{2}\til{\s}(\psi_t^2)~dB_t+\til{b}(\psi_t^2)~dt.
\end{equation}
Note that 0 is a solution of this S.D.E. (since $\til{\s}(0)=\til{b}(0)=0$). The solutions of this S.D.E. might be not unique.

\smallskip
Let $m(dx)$ be the speed measure of the diffusion~:
$$m(dx)=1_{]0,\pi[}(x)\exp\left[\int_{x_0}^x\frac{b(\p)}{\s^2(\p)}~d\p\right]~\frac{dx}{\s^2(x)}+m(\{0\})\d_0=g(x)~dx+m(\{0\})\d_0,$$
with $x_0\in ]0,\pi[$.

\begin{lem} If $\eta>\frac{1}{2}-\frac{(d-2)}{2\a}$, 0 is an entrance open boundary point.\end{lem}
\pf When $\eta>\frac{1}{2}-\frac{(d-2)}{2\a}$, 0 is an open boundary point. From Proposition 16.45 in \cite{br}, 0 is an entrance boundary point if and only if $\int_{0+}|s(x)|m(dx)<\infty$. For any positive $\e$, there exists a positive constant $D$ such that, for any $x\in]0,x_0[$,
$$|s(x)g(x)|\leq D~x^{(\frac{e}{c}-\a-\e)\wedge 0}x^{-\frac{e}{c}-\e+1}\leq D~x^{1-\a-2\e}.$$
This shows that $\int_{0+}s(x)m(dx)<\infty$ (choose $\e$ such that $2\e\leq 2-\a$). \qed

\medskip
This lemma implies that when $\eta>\frac{1}{2}-\frac{(d-2)}{2\a}$, there exists a positive $t$, a positive $\a$ and $p\in]0,1[$ such that for any $x\in]0,\pi[$, $P_x[\psi_t>\a]>p$. Proposition 6.8 implies that $S_t$ is not a flow of maps and since 0 is open, $S_t$ is diffusive without hitting.

\smallskip
Let now $d\in\{2,3\}$ (when $d\geq 4$ we always have $\eta>\frac{1}{2}-\frac{(d-2)}{2\a}$).

\begin{lem} If $\eta<\frac{1}{2}-\frac{(d-2)}{2\a}$, 0 is a closed boundary point. \end{lem}
\pf From Proposition 16.43 p.366 in \cite{br}, $T_0$ is finite or the boundary point 0 is closed if and only if for any $b\in ]0,\pi[$, $\int_0^b |s(x)-s(0)|m(dx)$ is finite.

We have $|s(x)-s(0)|g(x)\sim \int_0^x\exp\left[-\int_{x_0}^y\frac{b(\p)}{\s^2(\p)}~d\p\right]\times\frac{1}{\s^2(x)}\exp\left[\int_{x_0}^y\frac{b(\p)}{\s^2(\p)}~d\p\right]~dy$. Hence $|s(x)-s(0)|g(x)=O(x^{1-\a})$. This implies that $\int_0^b|s(x)-s(0)|m(dx)$ is finite. This proves that $T_0$ is finite a.s. \qed

\begin{lem} If $\eta<1-\frac{d}{\a^2}$, 0 is an exit boundary point. \end{lem}
\pf In \cite{br}, 0 is an exit boundary point if and only if $m(]0,x[)=\infty$ for all $x\in ]0,\pi[$. This is the case if $\mu-\a<-1$ (or if $\eta<1-\frac{d}{\a^2}$). Note that for $d=2$ or $3$ and $\a\in]0,2[$, $1-\frac{d}{\a^2}<\frac{1}{2}-\frac{(d-2)}{2\a}$. \qed

\medskip
Lemma 9.11 implies that when $\eta<1-\frac{d}{\a^2}$, the diffusion $\psi_t$ is absorbed at 0, and for any positive $r$, 
$$\lim_{d(x,y)\to 0}P^{(2)}_{(x,y)}[d(X_t,Y_t)>r]=\lim_{\p\to 0}P_{\p}[\psi_t>r]=0.$$
Now, applying proposition 6.7, we prove that the statistical solution is a flow of maps and this is a coalescent flow of maps (since 0 is an exit boundary point).

\begin{lem} If $\eta\in\left]1-\frac{d}{\a^2},\frac{1}{2}-\frac{(d-2)}{2\a}\right[$, 0 is a regular boundary point.\end{lem}
\pf In \cite{br}, we see that 0 is regular if  $m(]0,x[)<\infty$ for all $x\in ]0,\pi[$, which is the case when $\eta\in\left]1-\frac{d}{\a^2},\frac{1}{2}-\frac{(d-2)}{2\a}\right[$. \qed

\medskip
When $\eta\in\left]1-\frac{d}{\a^2},\frac{1}{2}-\frac{(d-2)}{2\a}\right[$, the two--point motion hits the diagonal. But there is no uniqueness of the solution of the S.D.E. satisfied by $\psi_t$ since 0 might be absorbing or (slowly or instantaneously) reflecting. In order to finish the proof of theorem 9.4, we are going to prove that 0 is instantaneously reflecting.

\smallskip
To prove this, for $\e\in ]0,1[$, let us introduce the covariance $C_\e=(1-\e)^2C$ (then, if $W_t$ is the cylindrical Brownian motion associated to $C$, $(1-\e)W_t$ is the cylindrical Brownian motion associated to $C_\e$) and $S_t^\e$ be the statistical solution associated to $P_t$ and $C_\e$.

For $f\in L^2(dx)$, $S_t^\e f=\sum_{n\geq 0}J^{n,\e}_tf$, where $J^{n,\e}_tf$ is the $n$th chaos in the chaos expansion of $S_t^\e f$ \footnote{Note that $S_t^\e=Q_{\log(1-\e)}S_t$, where $Q_\a$ is the Ornstein-Uhlenbeck operator on the Wiener space (used in Malliavin calculus~: see \cite{ma}).}.
 It is easy to see that $J^{n,\e}_tf=(1-\e)^nJ^{n}_tf$, where $J^{n}_tf$ is the $n$th chaos in the chaos expansion of $S_t f$, hence
\begin{equation}
E[(S_t^\e f-S_t f)^2]=\sum_{n\geq 1}(1-(1-\e)^{2n})E[(J^{n}_tf)^2].
\end{equation}
Hence it is clear that the $L^2(P)$-limit as $\e$ goes to 0 of $S_t^\e f$ is $S_tf$.

Let $(X_t^\e,Y_t^\e)$ be the Markov process associated to $P^{(2),\e}_t=E[S_t^{\e\,\otimes 2}]$ and $\psi_t^\e=d(X_t^\e,Y_t^\e)$. $\psi_t^\e$ is a diffusion with generator $L_\e$. It is easy to see that $L_\e=(1-(1-\e)^2)L_1+(1-\e)^2L$ (note that $A^{(2)}_\e=A\otimes I+I\otimes A+(1-\e)^2C=A^{(2)}_1+(1-\e)^2(A^{(2)}-A^{(2)}_1)$), and $L_\e=\s_{\e}^2(\p)\frac{d^2}{d\p^2}+b_{\e}(\p)\frac{d}{d\p}$, with
\begin{eqnarray}
\s_{\e}^2(\p)&=&(1-(1-\e)^2)\s_1^2(\p)+(1-\e)^2\s^2(\p),\\
b_{\e}(\p)&=&(1-(1-\e)^2)b_1(\p)+(1-\e)^2b(\p).
\end{eqnarray}
Let us remark that $L_1$ is the generator of the diffusion distance between two independent Brownian motions on $S^d$. Note that as $\p$ goes to 0,
\begin{equation} \s_1^2(\p)\sim\s_1^2(0)=2(a+b)KG(0)\qquad\hbox{and}\qquad b_1(\p)\sim\frac{2(d-1)}{\p}(a+b)KG(0)\end{equation}
and $\s_{\e}^2(\p)=(1-(1-\e)^2)\s_1^2(\p)(1+O(\p^\a))$ and $b_{\e}(\p)=(1-(1-\e)^2)b_1(\p)(1+O(\p^\a)$. Studying the scale function $s_\e$ of $\psi_t^\e$, we get that $s_\e(0+)=s_1(0+)=-\infty$ (as two independent Brownian motions cannot meet each other on $S^d$). We still have $s_\e(\pi-)=\infty$. Hence $\psi_t^\e\in ]0,\pi[$ for all positive $t$.

Let $m_\e$ be the speed measure of $\psi_t^\e$. Let $g_\e(x)=m_\e(dx)/dx$. As $m_\e(]0,\pi[)<\infty$, $m_\e$ is an invariant finite measure for the diffusion $\psi_t^\e$. As $\lim_{\e\to 0}\s^2_\e=\s^2$ and $\lim_{\e\to 0}b_\e=b$, we get that $\lim_{\e\to 0}g_\e(x)=g(x)$. Let us note $\e'=1-(1-\e)^2$ and let
\begin{equation}
f(\e',\p)=\frac{\e'b_1(\p)+(1-\e')b(\p)}{\e'\s^2_1(\p)+(1-\e')\s^2(\p)}.
\end{equation}
This function increases with $\e'$ if $\frac{b_1(\p)}{\s_1^2(\p)}\geq\frac{b(\p)}{\s^2(\p)}$. As $\frac{b_1(\p)}{\s_1^2(\p)}-\frac{b(\p)}{\s^2(\p)}\sim (d-1-\mu)\frac{1}{\p}$ as $\p$ goes to 0 and as $(d-1-\mu)$ is positive, there exists $\p_0$ such that for any $\p<\p_0$, $f(\e',\p)\geq\frac{b(\p)}{\s^2(\p)}=f(0,\p)$ and for $\e'<1/2$,
$$g_\e(x)\leq \frac{2}{\s^2(x)}\exp{\left(-\int_x^{\p_0}\frac{b(\p)}{\s^2(\p)}~d\p\right)}\times C_{\p_0},$$
where $C_{\p_0}=\sup_{\e\in [0,1]}\exp{(\int_{\p_0}^{x_0}f(\e',\p)~d\p)}<\infty$. The Lebesgue's dominated convergence theorem implies that $g_\e$ converges in $L^1([0,\pi])$ towards $g$.

Let $f$ and $g$ be continuous functions, then $E[f(X_t^\e)g(Y_t^\e)]=E[S_t^\e f(x)S_t^\e g(y)]$. Since $S_t^\e f$ and $S_t^\e g$ converge respectively towards $S_tf$ and $S_t g$ when $\e$ goes to 0 in $L^2(P)$, we get that $(X_t^\e,Y_t^\e)$ converges in distribution towards $(X_t,Y_t)$ when $\e$ goes to 0. This also implies that $\psi_t^\e$ converges in distribution towards $\psi_t$ when $\e$ goes to 0.

Since $m_\e$ is an invariant measure, for any continuous function $f$ on $[0,\pi]$, we have
\begin{equation}
\int E[f(\psi_t^\e)|\psi_0^\e=x]~m_\e(dx)=\int f~dm_\e.
\end{equation}
Since
\begin{eqnarray*}
\left|\int E[f(\psi_t^\e)|\psi_0^\e=x]~m_\e(dx)-\int E[f(\psi_t)|\psi_0=x]~g(x)dx\right|&\leq&\\
&&\hskip-302pt\leq~ \|f\|_\infty\int_0^\pi|g_\e(x)-g(x)|~dx+\left|\int_0^\pi \left(E[f(\psi_t^\e)|\psi_0^\e=x]-E[f(\psi_t)|\psi_0=x]\right)~g(x)dx\right|,
\end{eqnarray*}
we get that (because $g_\e$ converges in $L^1([0,\pi])$ towards $g$ and $\psi_t^\e$ converges in distribution toward $\psi_t$.)
$$\int E[f(\psi_t)|\psi_0=x]~m(dx)=\lim_{\e\to 0}\int E[f(\psi_t^\e)|\psi_0^\e=x]~m_\e(dx)=\lim_{\e\to 0}\int f~dm_\e=\int f~dm.$$
This implies that $g(x)dx$ is an invariant measure for $\psi_t$ and $m(dx)=g(x)dx$. Since $m(]0,x[)<\infty$ for all $x\in ]0,\pi[$, the diffusion $\psi_t$ is not absorbed in 0 and is reflected in 0.

\medskip
 In this case, 0 is a closed regular boundary point. This point is instantaneously reflecting since $m(\{0\})=0$. This implies the existence of a positive $t$, a positive $r$ and $p\in]0,1]$ such that for any $x\in ]0,\pi[$, $P_x[\psi_t\geq r]\geq p$. Then, applying proposition 6.8, the statistical solution is not a flow of maps. This finishes the proof of theorem 9.4. \qed

\bigskip
For $\a>2$, the statistical solution is an isotropic Brownian flow of diffeomorphisms. In Raimond \cite{ra}, the Lyapunov exponents of this flow are computed. The sign of the first Lyapunov exponent $\la_1(\a,d)$  describes the stability of the flow. It is unstable if $\la_1\geq 0$ and stable if $\la_1<0$. The computation of $\la_1(\a,d)$ gives
\begin{equation}
\la_1=\frac{(d-4)a+db}{d+2}\zeta(\a-1)+\left(\frac{d-1}{d+2}\right)[(d-4)a+db]\zeta(\a)-d\left(\frac{2(d-1)a+db}{d+2}\right)\zeta(\a+1),
\end{equation}
where $\zeta(\a)=\sum_{l\geq 1}\frac{1}{l^{\a}}$ is the zeta function. Therefore, we have $\la_1(\a,d)=0$ if and only if
\begin{equation}
\eta=\eta(\a,d)=\frac{-(d-4)\zeta(\a-1)-(d-1)(d-4)\zeta(\a)+2d(d-1)\zeta(\a+1)}{4\zeta(\a-1)+4(d-1)\zeta(\a)+d(d-2)\zeta(\a+1)}.
\end{equation}
It is easy to see that for fixed $\eta$, $\lim_{\a\to 2+}\la_1(\a,d)=+\infty$ if $d\geq 4$ or if $\eta>\frac{1}{2}-\frac{d-2}{4}=\frac{4-d}{4}$ and that $\lim_{\a\to 2+}\la_1(\a,d)=-\infty$ if $\eta<\frac{4-d}{4}$. Remark that $\lim_{\a\to 2-}1-\frac{d}{\a^2}=\lim_{\a\to 2-}\frac{1}{2}-\frac{(d-2)}{2\a}=\frac{4-d}{4}$. This shows that coalescence appears when $\la_1$ goes to $-\infty$ and splitting appears when $\la_1$ goes to $+\infty$.

\smallskip
The results of this section is given by phase diagrams in appendix B.

\setcounter{equation}{0}
\section{Isotropic statistical solution on $\R^d$.}

\subsection{Stationary and isotropic covariance functions on $\R^d$.}

On $\R^d$ with $d\geq 2$, the stationary isotropic covariance function $C$ are (see Le Jan \cite{lj}) such that $C^{ij}(x,y)=C^{ij}(x-y)$, for $(x,y)\in \R^d\times \R^d$, with
\begin{equation} C^{ij}(z)=\d^{ij}B_N(\|z\|)+\frac{z^iz^j}{\|z\|^2}(B_L(\|z\|)-B_N(\|z\|)), \end{equation}
with
\begin{eqnarray}
B_L(r)\!\!&=&\!\!\!\iint\cos(\rho u_1r)u_1^2\o(du)(F_L(d\rho)-F_N(d\rho))+\iint\cos(\rho u_1r)\o(du)F_N(d\rho),\\
B_N(r)\!\!&=&\!\!\!\iint\cos(\rho u_1r)u_2^2\o(du)(F_L(d\rho)-F_N(d\rho))+\iint\cos(\rho u_1r)\o(du)F_N(d\rho),
\end{eqnarray}
$F_L$ and $F_N$ being finite positive measures on $\R^+$. $\o(du)$ is the normalized Lebesgue measure on $S^{d-1}$. $F_L$ and $F_N$ represent respectively the gradient part and the zero divergence part of the associated Gaussian vector field.

\medskip
For $\a$ and $m$ positive reals, let $F(d\rho)=\frac{\rho^{d-1}}{(\rho^2+m^2)^{\frac{d+\a}{2}}}d\rho$, $F_L(d\rho)=aF(d\rho)$ and $F_N(d\rho)=\frac{b}{d-1}F(d\rho)$, where $a$ and $b$ are nonnegative. In the Fourier representation, (c is a positive constant)
\begin{equation}\hat{C}^{ij}(k)=c(\|k\|^2+m^2)^{-\frac{d+\a}{2}}\left(a\frac{k^ik^j}{\|k\|^2}+\frac{b}{d-1}\left(\d_{ij}-\frac{k^ik^j}{\|k\|^2}\right)\right).\end{equation}

\medskip
Notice that in the Fourier representation, the Laplace operator on vector fields is given by the multiplication by $-\|k\|^2$ and the projection $\pi$ on gradient vector fields (in the $L^2$ space) by $\frac{k^ik^j}{\|k\|^2}$. (i.e if $V$ is a vector field and $\hat{V}^i(k)$ its Fourier transform, $\hat{(\pi V)}^i(k)=\sum_j\frac{k^ik^j}{\|k\|^2}\hat{V}^j(k)$.)

\smallskip
Therefore, given a $L^2$ vector field, $U^j(y)=\int\sum_iC^{ij}(x-y)V^i(x)~dx$ can be expressed as $c(-\Delta+m^2)^{-\frac{d+\a}{2}}(a\pi V+\frac{b}{d-1}(I-\pi)V)$. Since $\l U,U\r_H=\l U,V\r_2=\int\l U(x),V(x)\r~dx$, the self--reproducing space appears to be the $L^2$-Sobolev space of order $s=\frac{d+\a}{2}$ (defined the same way as in secion 9.1) equipped with the norm
$$\|V\|^2=\frac{1}{a}\|\pi V\|^2_s+\frac{d-1}{b}\|(I-\pi)V\|^2_s,$$
where 
$$\|V\|^2_s=\frac{1}{c}\l (-\Delta+m^2)^sV,V\r_2.$$
Note that if $a$ or $b$ vanishes, the self--reproducing space is $H^{2,\frac{\a+d}{2}}$ restricted to divergence free vector fields or gradient vector fields.

\subsection{Phase transitions for the Sobolev statistical solution.}
Let $P_t$ be the semigroup of a Brownian motion on $\R^d$ with variance $(a+b)F(\R^+)$. Let $S_t$ be the statistical solution associated to $P_t$ and $C$. If $\a>2$, $C$ is $C^2$. Hence equation (8.2) is satisfied and the statistical solution $S_t$ is a flow of maps.

\smallskip
Suppose $a+b>0$ and let $\eta=\frac{b}{a+b}$. Then we have the theorem.

\begin{theo} For any $\a\in]0,2[$,
\begin{itemize}
\item For $d=2$ or $3$ and $\eta<1-\frac{d}{\a^2}$, the statistical solution is a coalescent flow of maps.
\item For $d=2$ or $3$ and $1-\frac{d}{\a^2}<\eta<\frac{1}{2}-\frac{(d-2)}{2\a}$, the statistical solution is diffusive with hitting.
\item For $d=2$ or $3$ and $\eta>\frac{1}{2}-\frac{(d-2)}{2\a}$ or for $d\geq 4$, the statistical solution is diffusive without hitting.
\end{itemize}
\end{theo}

\noindent{\bf Remark.} The results of this theorem are exactly the same as for the sphere.

\bigskip
\pf
Let us study the two--point motion $(X_t,Y_t)$ starting from $(x,y)$ (with $x\neq y$). Then $r_t=d(X_t,Y_t)$ is a diffusion in $\R^+$ (eventually reflected in 0), with generator $L=\s^2(r)\frac{d^2}{dr^2}+b(r)\frac{d}{dr}$ (see Le Jan \cite{lj}), with
\begin{eqnarray}
\s^2(r)&=&B-B_L(r),\\
b(r)&=&(d-1)\frac{B-B_N(r)}{r},
\end{eqnarray}
where $B=B_L(0)=B_N(0)=\frac{a+b}{d}F(\R^+)$.

\begin{lem} For $\a\in ]0,2[$, as $r$ goes to 0,
\begin{description}
\item{i)} $\iint\cos(\rho u_1r)\o(du)F(d\rho)=F(\R^+)-\a_1r^\a+o(r^\a)$.
\item{ii)} $\iint\cos(\rho u_1r)u_1^2\o(du)F(d\rho)=\frac{F(\R^+)}{d}-\a_2r^\a+o(r^\a)$.
\item{iii)} $\iint\cos(\rho u_1r)u_2^2\o(du)F(d\rho)=\frac{F(\R^+)}{d}-\a_3r^\a+o(r^\a)$.
\end{description}
with $\a_2=\frac{\a+1}{d+\a}\a_1$, $\a_3=\frac{1}{d+\a}\a_1$ and
$$\a_1=c_d\left(\int_0^\infty(1-\cos x)\frac{dx}{x^{\a+1}}\right)\left(\int_0^{\frac{\pi}{2}}(\cos\t)^\a(\sin\t)^{d-2}d\t\right).$$
\end{lem}

\pf 
For $r>0$, making the change of variable $x=\rho u_1r$,
\begin{eqnarray*}
\iint(1-\cos(\rho u_1r))\o(du)F(d\rho) 
&=& c_d\int_0^1\int_0^\infty(1-\cos(\rho u_1r))(1-u_1^2)^{\frac{d-2}{2}}du_1\frac{\rho^{d-1}d\rho}{(\rho^2+m^2)^{\frac{d+\a}{2}}}\\
=~ r^\a c_d\hskip-10pt&&\hskip-17pt\int_0^1\left(\int_0^\infty(1-\cos x)\frac{x^{d-1}dx}{(x^2+r^2u_1^2m^2)^{\frac{d+\a}{2}}}\right)u_1^\a(1-u_1^2)^{\frac{d-2}{2}}du_1.
\end{eqnarray*}
As $\lim_{r\to 0}\int_0^\infty(1-\cos x)\frac{x^{d-1}}{(x^2+r^2u_1^2m^2)^{\frac{d+\a}{2}}}dx=\int_0^\infty(1-\cos x)\frac{dx}{x^{\a+1}}<\infty$, we get that 
$$\lim_{r\to 0}\frac{1}{r^\a}\iint(1-\cos(\rho u_1r))\o(du)F(d\rho)=c_d\left(\int_0^\infty(1-\cos x)\frac{dx}{x^{\a+1}}\right)I(d-2,\a)=\a_1,$$
with $I(n,t)=\int_0^{\frac{\pi}{2}}(\cos\t)^t(\sin\t)^n~d\t=\frac{1}{2}B(\frac{n+1}{2},\frac{t+1}{2})$ for $t\geq 0$ and $n\in\N$, and $B(x,y)=\frac{\G(x)\G(y)}{\G(x+y)}$. This shows i). ii) and iii) can be obtained the same way with
$$\a_2=c_d\int_0^\infty(1-\cos x)\frac{dx}{x^{\a+1}}I(d-2,\a+2)$$
and $\a_1=\a_2+(d-1)\a_3$ (note that $\int u_1^2\o(du)=\frac{1}{d}$). It is easy to see that for $\a>0$ and $d\geq 1$,
$$I(d-2,\a+2)=\frac{\a+1}{d+\a}I(d-2,\a).$$
Therefore, $\a_2=\frac{\a+1}{d+\a}\a_1$. With the relation $\a_1=\a_2+(d-1)\a_3$, we get that $\a_3=\frac{1}{d+\a}\a_1$. \qed

\begin{rk} As $z$ goes to 0,
$$C^{ij}(z)=B\d^{ij}-\frac{\a_1}{d-1}\left[((d-1)a+(d+\a-1)b)\d^{ij}-\a((d-1)a-b)\frac{z^iz^j}{\|z\|^2}\right]\|z\|^\a(1+o(1)),$$
Let us note that the dependence on $m$ only appears in $B$.
\end{rk}

From this lemma, it is easy to see that as $r$ goes to 0,
\begin{eqnarray}
\s^2(r)&=&\frac{(a+b)\a_1}{d+\a}(\a+1-\a\eta)r^\a(1+o(1)),\\
b(r)&=&\frac{(a+b)\a_1}{d+\a}(d-1+\a\eta)r^{\a-1}(1+o(1)).
\end{eqnarray}
Note that we get the same behaviour of $\s$ and $b$ around 0 as in section 9.2.

\smallskip
As in section 9.2, let us study $s$, the scale function of the diffusion $r_t$.

Since $B_L(r)$ and $B_N(r)$ converge towards 0 as $r$ goes to $\infty$ (as Fourier transforms of finite measures), we get that as $r$ goes to $\infty$, $\log(s'(r))\sim (1-d)\log(r)$. Therefore $s(+\infty)$ is finite if and only if $d\geq 3$.

We also see that $s(0+)=-\infty$ if $\eta>\frac{1}{2}-\frac{(d-2)}{2\a}$ and $s(0+)$ is finite if $\eta<\frac{1}{2}-\frac{(d-2)}{2\a}$. 

\smallskip
Let $m$ be the speed measure of the diffusion. 
Let us study the boundary point 0. 

As $m(]0,x[)<-\infty$ for any positive $x$ if $\eta>1-\frac{d}{\a^2}$, as in section 9.2, with a similar proof, we can prove that if $\eta\in ]1-\frac{d}{\a^2},\frac{1}{2}-\frac{(d-2)}{2\a}[$, the diffusion $r_t$ is instantaneously reflecting at 0. The only thing there is to change in the proof is to take the test function $f$ in (9.28) with compact support and to remark that $g_\e$ converges towards $g$ in $L^1_{loc}(\R^+)$.

If $\eta<1-\frac{d}{\a^2}$ (note that $1-\frac{d}{\a^2}\leq\frac{1}{2}-\frac{d-2}{\a}$), 0 is an exit boundary point and the diffusion is absorbed by 0.

\smallskip
Therefore, we get that 
\begin{itemize}
\item If $d\geq 3$ and $\eta\in ]1-\frac{d}{\a^2},\frac{1}{2}-\frac{(d-2)}{2\a}[$, $r_t$ is instantaneously reflecting at 0 and is transient. In this case, as in section 9.2, $(S_t)_{t\geq 0}$ is diffusive with hitting.
\item If $d=2$ and $\eta\in ]1-\frac{d}{\a^2},\frac{1}{2}-\frac{(d-2)}{2\a}[$, $r_t$ is instantaneously reflecting at 0 and is recurrent. In this case, as in section 9.2, $(S_t)_{t\geq 0}$ is diffusive with hitting.
\item If $d\geq 3$ and $\eta<1-\frac{d}{\a^2}$, $r_t$ is absorbed at 0 with probability $\frac{s(\infty)-s(r_0)}{s(\infty)-s(0)}$ and converges towards $+\infty$ with probability $\frac{s(r_0)-s(0)}{s(\infty)-s(0)}$. In this case, as in the section 9.2, $(S_t)_{t\geq 0}$ is a coalescent flow of maps.
\item If $d=2$ and $\eta<1-\frac{d}{\a^2}$, $r_t$ is absorbed at 0 a.s. In this case, as in section 9.2, $(S_t)_{t\geq 0}$ is a coalescent flow of maps.
\end{itemize}

If $\eta>\frac{1}{2}-\frac{d-2}{2\a}$, then we have that $s(0)=-\infty$. In this case, 0 is an entrance boundary point as $\int_{0+}|s(x)|dm(x)<\infty$. $r_t$ is recurrent if $d= 2$ and transient if $d\geq 3$. As in the section 9.2, we prove that $(S_t)_{t\geq 0}$ is diffusive without hitting. \qed

\medskip
For $\a>2$, the statistical solution is a stationary isotropic Brownian flow of diffeomorphisms. In Le Jan \cite{lj}, the Lyapunov exponents of this flow are computed. The sign of the first Lyapunov exponent $\la_1(\a,d)$  describes the stability of the flow. It is unstable if $\la_1\geq 0$ and stable if $\la_1<0$. The computation of $\la_1(\a,d)$ gives (see \cite{lj})
\begin{equation}
\la_1=\frac{1}{2(d+2)}((d-4)a+db)\int \rho^2 F(d\rho),
\end{equation}
Therefore, we have $\la_1(\a,d)=0$ if and only if $d\leq 4$ and
\begin{equation}
\eta=\eta(d)=\frac{4-d}{4}
\end{equation}
As in the section 9.2, we see that for fixed $\eta$, $\lim_{\a\to 2+}\la_1(\a,d)=+\infty$ if $d\geq 4$ or if $\eta>\frac{1}{2}-\frac{d-2}{4}=\frac{4-d}{4}$ and that $\lim_{\a\to 2+}\la_1(\a,d)=-\infty$ if $\eta<\frac{4-d}{4}$. This shows that coalescence appears when $\la_1$ goes to $-\infty$ and splitting appears when $\la_1$ goes to $+\infty$.

\smallskip
Remark that $\lim_{\a\to 2-}1-\frac{d}{\a^2}=\lim_{\a\to 2-}\frac{1}{2}-\frac{(d-2)}{2\a}=\frac{4-d}{4}$.

\smallskip
The results of this section are given by phase diagrams in appendix B.

\setcounter{equation}{0}
\section{Reflecting flows.}

Let $D$ be an open convex domain in $\R^d$ with $C^1$ boundary $\pa D$. Let $d$ be the Euclidean metric in $\R^d$. For any $x\in\pa D$, we denote $n(x)$ the directed inward unit normal vector to $\pa D$.

Let $P_t$ be the semigroup of the Brownian motion in $D$ reflected on $\partial D$. $P_t$ is associated to the Dirichlet form $(\E,\F)$, where $\F=H^1(D)=\{f\in L^2(D,dx),~|\nabla f|\in L^2(D,dx)\}$ equipped with the form $\frac{1}{2}\int_D |\nabla f|^2~dx$ (see \cite{fuku}, 1.3.2). Let $C(x,y)$ be a covariance function in $D\times D$ such that $C^{ij}(x,x)=\delta^{ij}$ and satisfying (8.1).

We can construct a statistical solution associated to $P_t$ and $C$. Let $P^{(2)}_t$ be the semigroup of the two--point motion  $(X_t,Y_t)$. Let $P^{(2)}_{(.,.)}$ be the law of the two--point motion.

 We know that $X_t$ and $Y_t$ are two diffusions in $D$ reflected on $\pa D$. Let $\p_t$ and $\psi_t$ denote the local times of $X_t$ and $Y_t$ on $\pa D$.

\begin{lem}For $h(x,y)=d^2(x,y)$,
$P^{(2)}_t h(x,y)\leq h(x,y)e^{Ct}$.
\end{lem}

\pf  Let us note 
$$L^{(2)}=A_x+A_y+\sum_{i,j}C^{ij}(x,y)\pa_{x_i}\pa_{x_j}.$$
From (8.1) and the Lipschitz conditions, we get that
$$L^{(2)}h(x,y)\leq C~h(x,y).$$
Using Tanaka's formula, there exists a martingale $M_t$ such that
\begin{eqnarray}
h(X_t,Y_t)-h(x,y)
&=& M_t+\int_0^tL^{(2)}h(X_s,Y_s)~ds\\
&+&\int_0^t\l\n_x h(X_s,Y_s),n(X_s)\r d\p_s+\int_0^t\l\n_y h(X_s,Y_s),n(Y_s)\r d\psi_s.
\end{eqnarray}
As $\n_x h(x,y)=2(x-y)$, using the fact that $D$ is convex, we get that for $x\in\pa D$
$$\l\n_x h(x,y),n(x)\r<0.$$
This implies that
$$h(X_t,Y_t)-h(x,y)\leq M_t+C\int_0^t h(X_s,Y_s)~ds.$$
Taking the expectation, we get that $P^{(2)}_t h(x,y)-h(x,y)\leq C\int_0^t P^{(2)}_s h(x,y)~ds$. Hence the lemma. \qed

\begin{theo}
The statistical solution is a flow of maps.
\end{theo}

\pf This is the same proof as the proof of theorem 8.2. \qed

\appendix
\renewcommand {\theequation}{A.\arabic{equation}}
\setcounter{equation}{0}
\section{Proof of lemma 9.5.}
Take $\p\in]0,\pi[$. At first, we are going to prove that $I(\p)=\sum_{l\geq 2}\frac{1}{(l-1)^{\a+1}}\left|\frac{d}{d\p}\g_l(\cos\p)\right|$ is finite. As $\frac{1}{l^{\a}}=\int_0^\infty e^{-ls}s^{\a-1}~\frac{ds}{\G(\a)}$,
\begin{eqnarray}
I(\p)&\leq&\int_0^\pi\int_0^\infty\sum_{l\geq 1}[e^{-s}|z(\p,\t)|]^l\frac{\left|\frac{d}{d\p}z(\p,\t)\right|}{|z(\p,\t)|}s^{\a-1}\sin^d\t~\frac{ds~d\t}{\G(\a)c_d}\\
&\leq&\int_0^\pi \int_0^\infty f_{\p,\t}(s)~ds~d\t=2\int_0^{\frac{\pi}{2}} \int_0^\infty f_{\p,\t}(s)~ds~d\t,
\end{eqnarray}
with $f_{\p,\t}(s)=\frac{e^{-s}\left|\frac{d}{d\p}z(\p,\t)\right|}{1-e^{-s}|z(\p,\t)|}\frac{s^{\a-1}}{\G(\a)c_d}$. It is easy to see that
\begin{equation}
\int_1^\infty f_{\p,\t}(s)~ds\leq\frac{1}{\G(\a)c_d}\int_1^\infty\frac{e^{-s}s^{\a-1}}{(1-e^{-s})}ds<\infty.\end{equation}
On the other hand,
\begin{equation}
\int_0^1 f_{\p,\t}(s)~ds\leq \frac{1}{\G(\a)c_d}\int_0^1\frac{ds}{1-e^{-s}|z(\p,\t)|}=\frac{1}{\G(\a)c_d}F_\p(\t).
\end{equation}
Let $x_\p(\t)=-\log|z(\p,\t)|$, then $F_\p(\t)=\int_{x_\p(\t)}^{x_\p(\t)+1}\frac{dt}{1-e^{-t}}$. As $\lim_{\t\to 0+}x_\p(\t)=0$, we have $F_\p(\t)\sim -\log x_\p(\t)$ as $\t$ goes to 0. From this, we see that $F_\p(\t)=O(\log\t)$ as $\t$ goes to 0. This implies that $I(\p)$ is finite. 

Now, applying the derivation under the integral theorem, we prove that $G$ is differentiable on $]0,\pi[$ and that for $\p\in]0,\pi[$,
\begin{eqnarray}G'(\p)&=&\sum_{l\geq 1}\int_0^\pi\frac{[z(\p,\t)]^{l-1}\frac{d}{d\p}z(\p,\t)}{l^{\a}}\sin^d\t~\frac{d\t}{c_d}\\
&=&\int_0^\pi\int_0^\infty\sum_{l\geq 1}[e^{-s}z(\p,\t)]^l\frac{\frac{d}{d\p}z(\p,\t)}{z(\p,\t)}s^{\a-1}\sin^d\t~\frac{ds~d\t}{\G(\a)c_d}\\
&=&\int_0^\pi\int_0^\infty \frac{e^{-s}\frac{d}{d\p}z(\p,\t)}{1-e^{-s}z(\p,\t)}s^{\a-1}\sin^d\t~\frac{ds~d\t}{\G(\a)c_d}.
\end{eqnarray}
As $z(\p,\pi-\t)=\overline{z(\p,\t)}$,
$$G'(\p)=-\int_0^{\frac{\pi}{2}}\int_0^\infty\frac{a(s,\p)-\sin^2\t}{b(s,\p)+\cos^2\t}\times\frac{\cos\p}{\sin\p}s^{\a-1}\sin^d\t~\frac{2ds~d\t}{\G(\a)c_d},$$
with $a(s,\p)=\frac{1}{e^{-s}\cos\p}$ and $b(s,\p)=\frac{(1-e^{-s}\cos\p)^2}{e^{-2s}\sin^2\p}$. Changing of variables ($s=t\p$),
\begin{eqnarray}
-\frac{G'(\p)}{\p^{\a-1}} &=& \int_0^{\frac{\pi}{2}}\int_0^\infty \frac{a(t\p,\p)-\sin^2\t}{b(t\p,\p)+\cos^2\t}\times\frac{\p\cos\p}{\sin\p}t^{\a-1}\sin^d\t~\frac{2dt~d\t}{\G(\a)c_d}\\
&=& \int_0^{\frac{\pi}{2}}\int_0^\infty I(t,\p,\t)~dt~d\t.
\end{eqnarray}
Let $\e>0$, there exists a positive constant $C_\e$ such that for any $t\in[0,\e]$,
\begin{equation} 0\leq I(t,\p,\t)\leq C_\e t^{\a-1}. \end{equation}

Remark also that
\begin{equation} I(t,\p,\t)\leq C_{d,\a}\frac{t^2\p^2e^{-t\p}}{(1-e^{-t\p})^2}\times t^{\a-3},\end{equation}
where $C_{d,\a}$ is a positive constant. Let $C=C_{d,\a}\sup_{x>0}\frac{x^2e^{-x}}{(1-e^{-x})^2}<\infty$, then, for any positive $t$
\begin{equation} 0\leq I(t,\p,\t)\leq C~t^{\a-3}.\end{equation}

As $F(t)=C_\e t^{\a-1}1_{0<t\leq\e}+C~t^{\a-3}1_{t>\e}$ belongs to $L^1(d\t\otimes dt)$ for $\a\in]0,2[$, $\lim_{\p\to 0}a(t\p,\p)=1$ and $\lim_{\p\to 0}b(t\p,\p)=t^2$, by the Lebesgue dominated convergence theorem,
\begin{equation} \lim_{\p\to 0}\frac{G'(\p)}{\p^{\a-1}}=-\int_0^{\frac{\pi}{2}}\int_0^\infty\frac{\cos^2\t}{t^2+\cos^2\t}t^{\a-1}\sin^d\t~\frac{2d\t~dt}{c_d\Gamma(\a)}=-\a K.\end{equation}
We have proved the second limit. The first limit is easy to obtain as
\begin{eqnarray*}
G(0)-G(\p) &=& -\int_0^\p G'(x)~dx\\
&=& -K~\p^{\a}+o(\p^\a).
\end{eqnarray*}
This finishes the proof of the lemma. \qed

\newpage

\section{Phase diagrams for the Sobolev statistical solutions.}

\bigskip
\centerline{\epsfysize=10cm\epsfbox{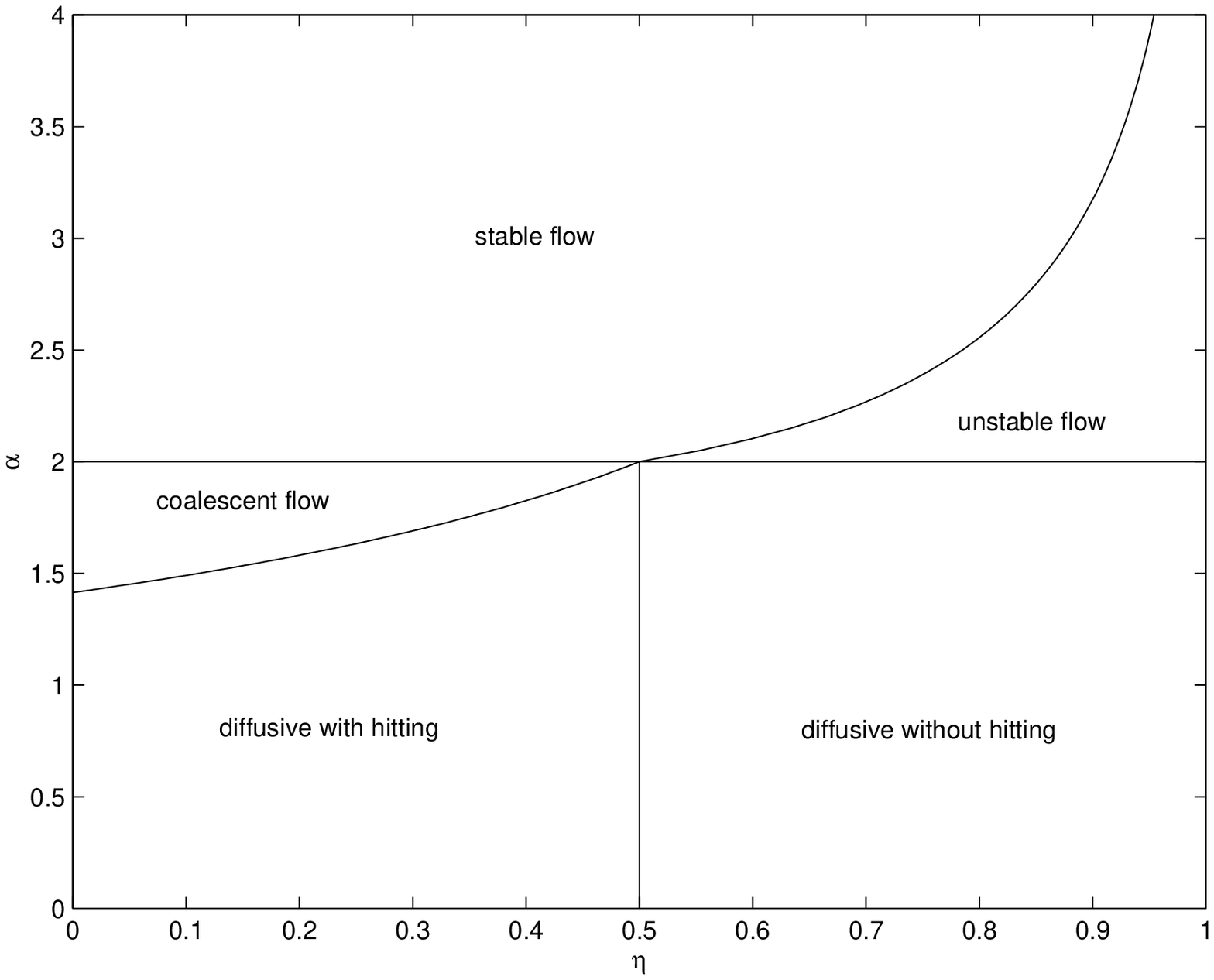}}
\centerline{{\bf Figure 1~:} Phase diagram on $S^2$.}

\bigskip
\centerline{\epsfysize=10cm\epsfbox{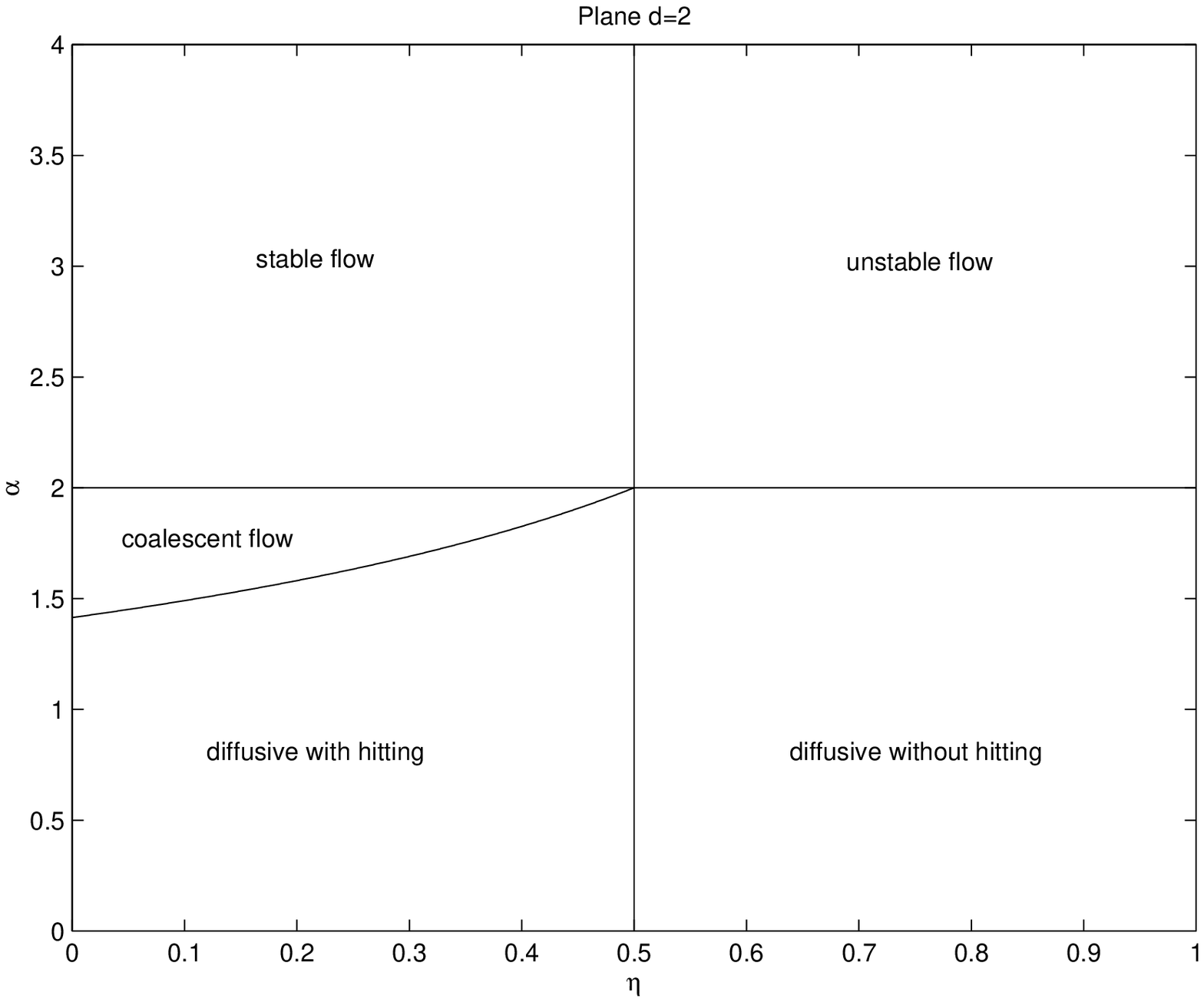}}
\centerline{{\bf Figure 2~:} Phase diagram on $\R^2$.}

\bigskip
\centerline{\epsfysize=10cm\epsfbox{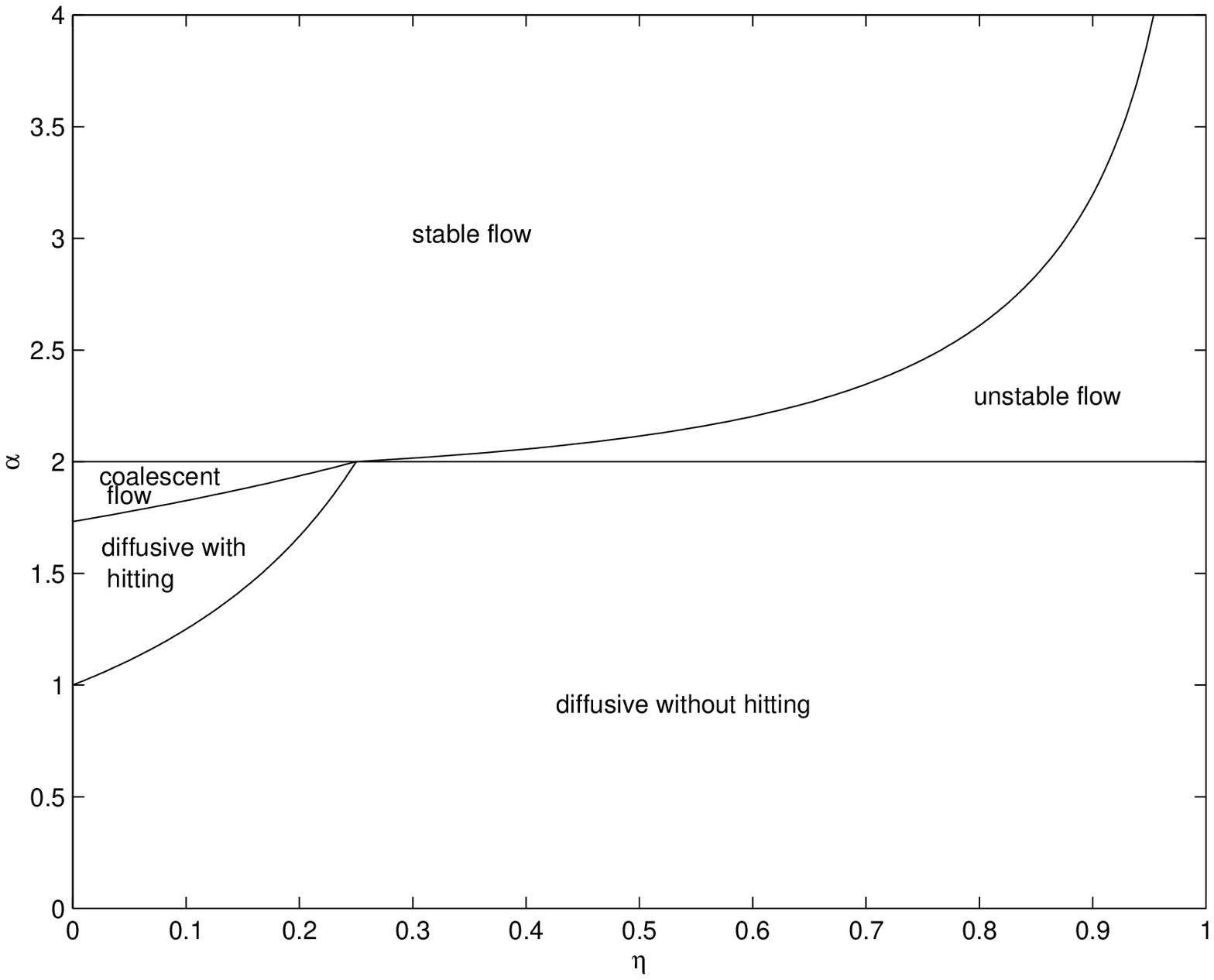}}
\centerline{{\bf Figure 3~:} Phase diagram on $S^3$.}

\bigskip
\centerline{\epsfysize=10cm\epsfbox{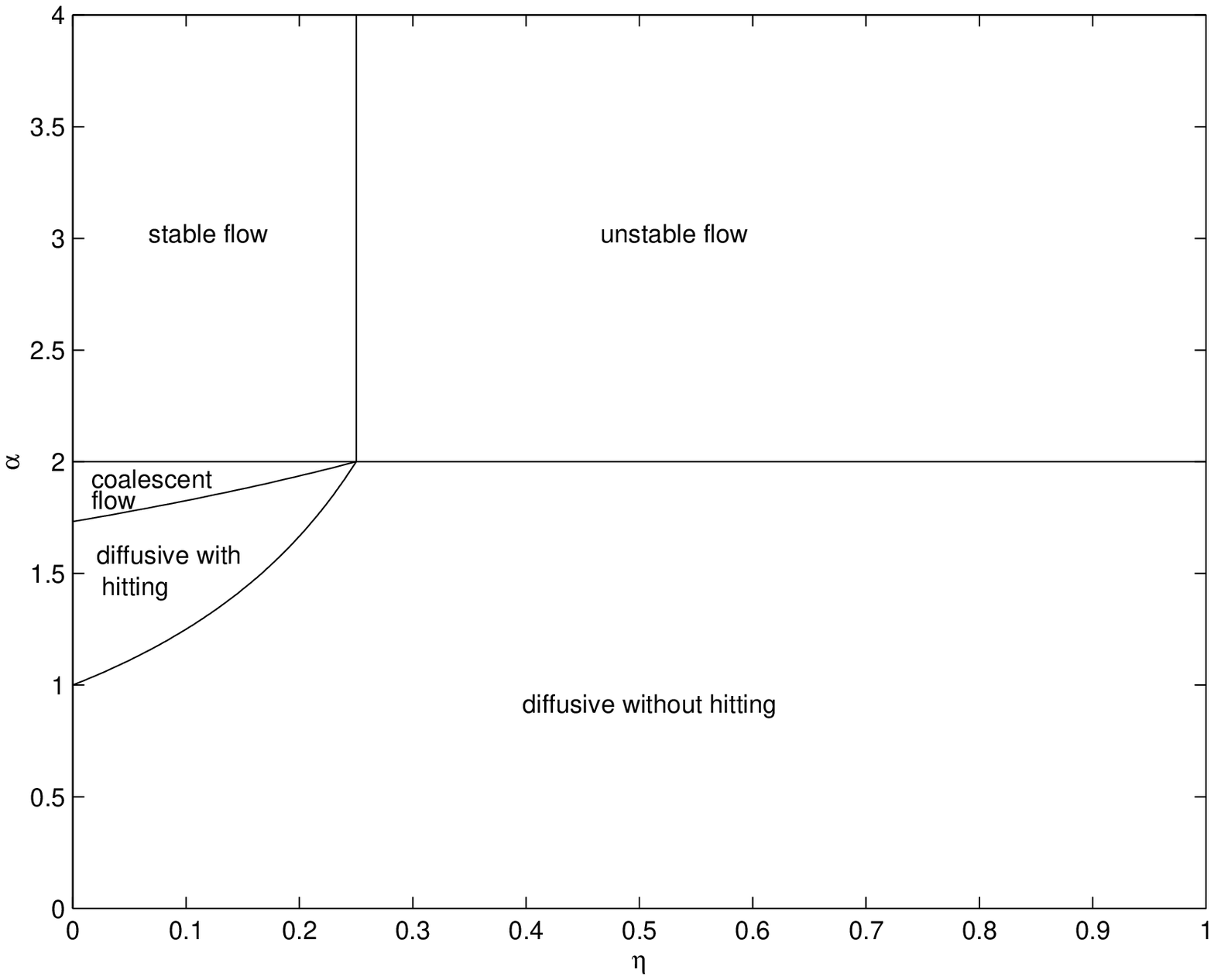}}
\centerline{{\bf Figure 4~:} Phase diagram on $\R^3$.}

\bigskip
\centerline{\epsfysize=10cm\epsfbox{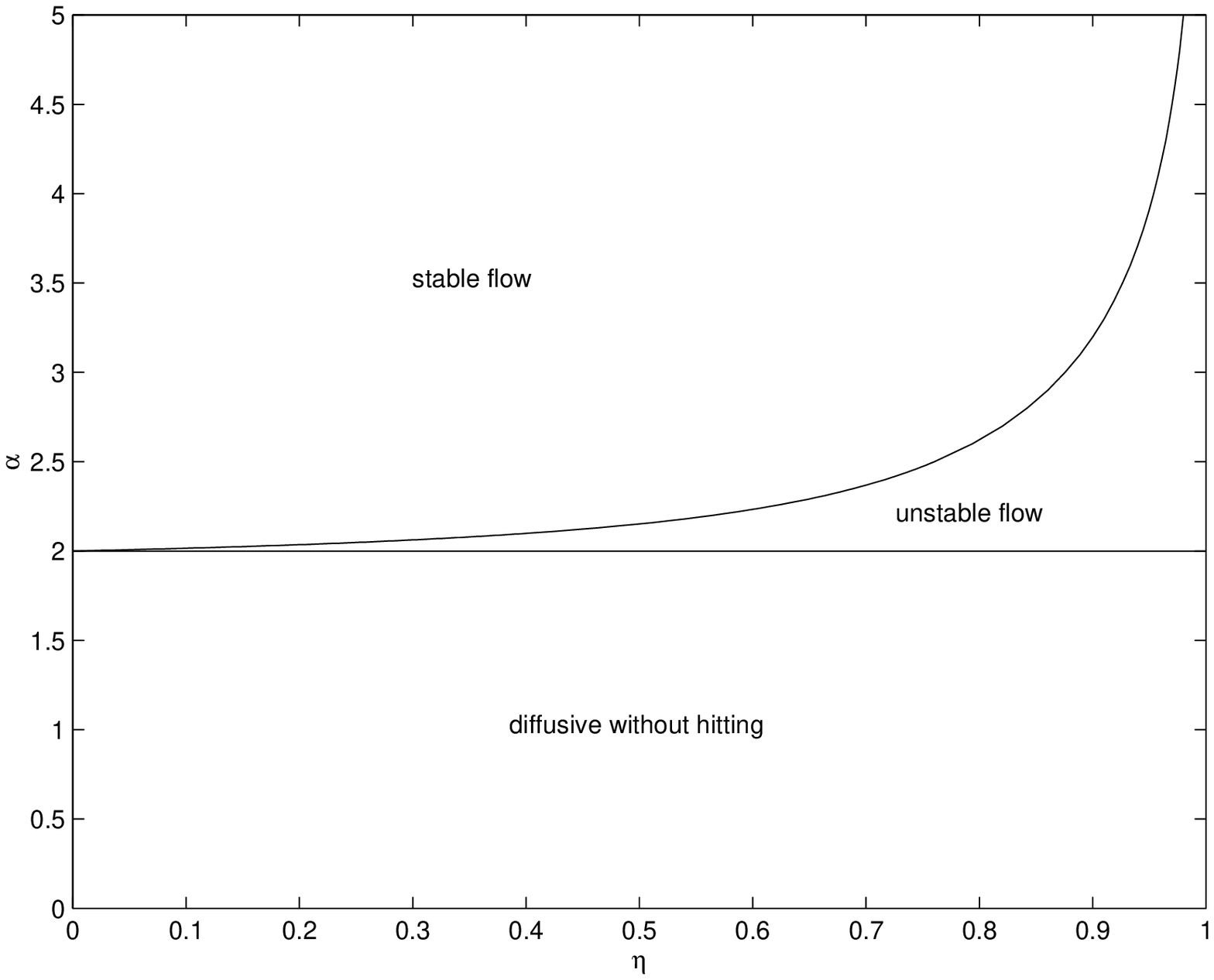}}
\centerline{{\bf Figure 5~:} Phase diagram on $S^4$.}

\bigskip
\centerline{\epsfysize=10cm\epsfbox{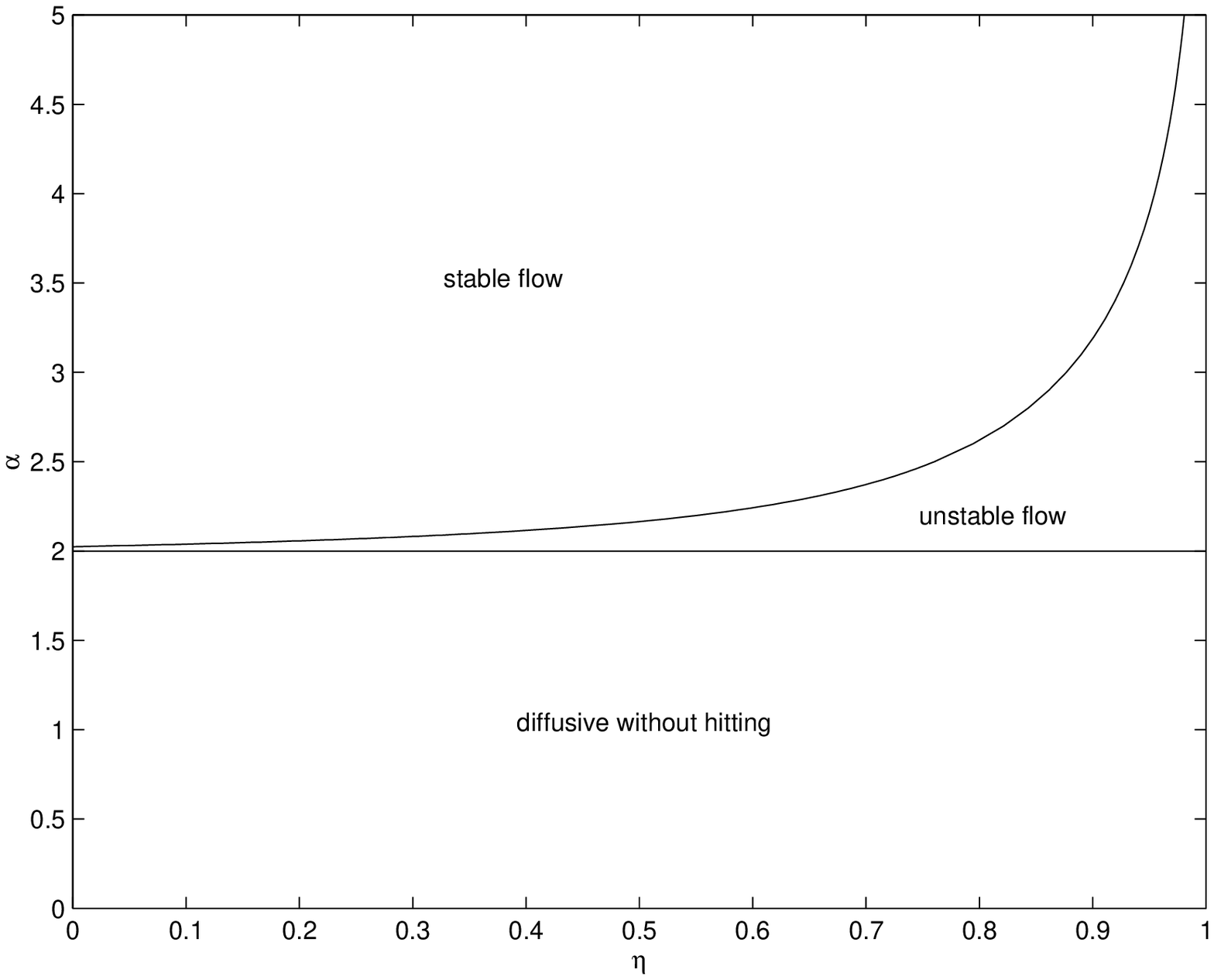}}
\centerline{{\bf Figure 6~:} Phase diagram on $S^5$.}

\bigskip
\centerline{\epsfysize=10cm\epsfbox{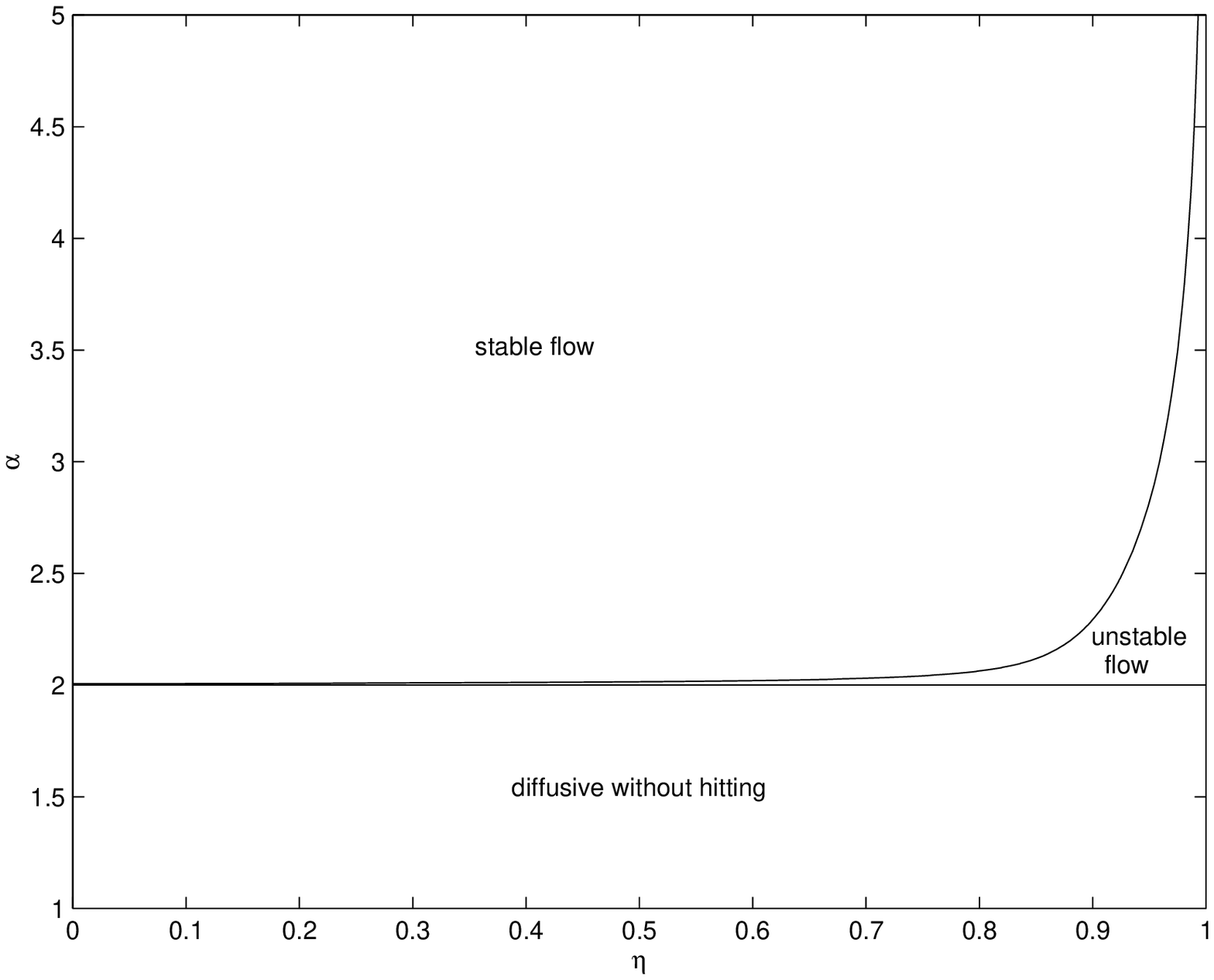}}
\centerline{{\bf Figure 7~:} Phase diagram on $S^{50}$.}

\bigskip
Let us remark that when $\a<2$, the diagrams are exactely the same for the sphere and for the plane. For the sphere, we see that, for $\a>2$ and $\eta\leq 2-\frac{\zeta(\a)}{\zeta(\a+1)}$, the flow gets stable when $d$ goes to $\infty$~: (9.30) implies that $\lim_{d\to\infty}\eta(\a,d)=2-\frac{\zeta(\a)}{\zeta(\a+1)}$ for $\a>2$. We see that, for any $d$ and $\eta\in [0,1[$, the flow gets stable when $\a$ goes to $\infty$~: (9.30) implies that $\lim_{\a\to\infty}\eta(\a,d)=1$.

\end{document}